\documentclass[12pt,a4paper,twoside,final,notitlepage, leqno]{article}
\usepackage[english]{babel}
\usepackage[T1]{fontenc}  
\usepackage{graphicx}
\usepackage{amsmath,amsthm,epsfig,amsfonts,bbm}  

\newcommand{\pequationdeb}{$$ \left\{ \begin{minipage}[c]{130mm}}
\newcommand{\pequationfin}{\end{minipage}
                           \right. $$}

\newcommand{\moneq}{\vspace*{-15pt} \begin{equation} \displaystyle } 
\newcommand{\nototo}{ \vspace*{-6pt} \noindent }
 
\newcommand{\beq}     {\begin{equation}}
\newcommand{\enq}     {\end{equation}}
\newcommand{\be}    {\begin{enumerate}}
\newcommand{\ee}    {\end{enumerate}}

\newcommand{\Bb}

\def\resume{\if@twocolumn
\section*{R\'esum\'e}
\else \small
\quotation{\bf \it R\'esum\'e \rule[1mm]{1.5mm}{0.2mm}\vspace{0pt}}
\fi}
\def\endresume{\if@twocolumn\else\endquotation\fi}
\def\abstract{\if@twocolumn
\noindent\section*{{\bf Abstract}}
\else \small
\quotation{\noindent \bf {Abstract.} \rule[1mm]{1.5mm}{0.2mm}\vspace{0pt}}
\fi}
\def\endabstract{\if@twocolumn\else\endquotation\fi}

\newtheorem{theorem}{Theorem}
\newtheorem{proposition}[theorem]{Proposition}

\def\section*#1{}

\begin{document}  
\title{ \bf \LARGE     ~ \vspace{1cm}  ~\\    Lorentz Transform    ~\\  
and Staggered Finite Differences     ~\\  
  for Advective Acoustics   ~\\~  } 

\author { { \large  F. Dubois~$^{ab}$, E. Duceau~$^{c}$,    F. Mar\'echal~$^{cd}$ 
and I. Terrasse~$^{c}$}     \\ ~\\
{\it \small  $^a$  Conservatoire National des Arts et M\'etiers, Saint Cyr l'Ecole,
  France.} \\
{\it \small  $^b$  Applications Scientifiques du Calcul Intensif, Orsay, France.}\\
{\it \small  $^c$   European Aeronautics Defence and Space Compagny,}\\ 
{\it \small  Research Center,  Suresnes, France.} \\  
{\it \small  $^d$   Ecole Nationale des Ponts et Chauss\'ees, Marne-la-Vall\'ee, France. }}   
\date{ Ao\^ut  2002~\protect\footnote{~Edition  06 may 2011.}}

\maketitle
 
\begin{abstract}
We study acoustic wave propagation in a uniform stationary flow. We develop a
  method founded on the Lorentz transform and a hypothesis of irrotationality of the
  acoustic perturbation. After a transformation of the space-time and of the unknown
  fields, we derive a system of partial differential equations that eliminates the
  external flow and deals with the classical case of non advective acoustics. A sequel of
  the analysis is a new set of perfectly matched layers equations in the spirit of the
  work of Berenger and Collino. The numerical implementation of the previous ideas is
  presented with the finite differences method HaWAY on cartesian staggered
  grids. Relevant numerical tests are proposed.
$ $ \\[5mm]
 {\it  Keywords}:  Perfectly matched layers, finite differences, HaWAY method.
$ $ \\[5mm]
 {\it  AMS classification}: 65M06,  76N15
\end{abstract}

\newpage 
 
\bigskip \bigskip  \noindent {\bf \large 1) \quad Introduction}

\noindent 
The acoustic dimensioning of a civil aircraft requires the use of numerical 
models to predict the radiated acoustic field emitted by the engine. 
We study this problem in the case of advective acoustics. The specific case of 
classical acoustic wave propagation can be viewed as a scalar version of Maxwell 
equations. We use our previous experience acquired in the simulation of the 
propagation of electromagnetic waves using the finite differences method \cite{cray} 
to rapidly develop a three-dimensional software simulating the propagation of acoustic waves.

\smallskip \noindent 
Our work is structured as follows. In order to take into account the external 
aerodynamic flow, we first come back to the equations of gas dynamics and consider 
the acoustic field as a first order linear perturbation of such a flow. Then we use 
physical ideas based on the Lorentz group invariance in Section 3,  
in the spirit of \cite{astley}, to deal in Section 4 
with the case of advective 
acoustics in the same way as the non-advective ones. We develop in Section 5 
a smart solution of the difficult problem of absorbing layers. The numerical aspect 
with the use of staggered grids is tackled in Section 6,  
and relevant 
physical and numerical tests are proposed in Section 7 

\bigskip \bigskip  \noindent {\bf \large 2) \quad Non linear acoustics} 

  \bigskip  \noindent {\bf \large 2-1   \quad  Barotrope gas dynamics}

\noindent 
We consider the propagation of sound waves in a uniform two-dimensional subsonic 
flow of a compressible fluid. This phenomenon is described by the nonlinear Euler 
equations for gas dynamics, see \cite{landau} for example, as~:

\moneq
\left \lbrace
  \begin{array}{l}
    \displaystyle{ \frac{\partial \check{\rho}}{\partial t}} + \mbox{div}~({\check{\rho} 
\check{\textbf{u}}}) = 0 \\
    \displaystyle{ \frac{\partial (\check{\rho} \check{u})}{\partial t} + 
\frac{\partial}{\partial x} (\check{\rho} \check{u}^2 +\check{p}) + 
\frac{\partial}{\partial y} (\check{\rho} \check{u} \check{v}) = 0 } \\
    \displaystyle{ \frac{\partial (\check{\rho} \check{v})}{\partial t} + 
\frac{\partial}{\partial x} (\check{\rho} \check{u} \check{v}) + 
\frac{\partial}{\partial y} (\check{\rho} \check{v}^2 +\check{p}) = 0 }  \\
    \displaystyle{ \frac{\partial \check{s}}{\partial t} + \check{u} 
\frac{\partial \check{s}}{\partial x} + \check{v} \frac{\partial \check{s}}{\partial y} = 0 }~,
  \end{array}
\right.
\label{sys1_1} 
\end{equation}

\noindent where $\check{{\textbf{u}}}=(\check{u}, \check{v})$ is the velocity 
vector, $\check{\rho}$ the density of the fluid, $\check{p}$ the pressure of 
the fluid and $\check{s}$ the entropy. We also know that~:

\moneq
        \displaystyle{ \frac{\check{p}}{\check{p_0}} = 
\frac{\check{\rho}^{\gamma}}{\check{\rho_0}^{\gamma}} \exp( \frac{\check{s}}{C_{V}} ) }~,
\label{sys1_2}
\end{equation}

\noindent where $C_{V}$ is the calorific capacity at constant volume 
and $(\rho_0, p_0)$ a state of reference.

\bigskip  \newpage   \noindent {\bf \large 2-2   \quad  Linearization around a stationary state}

\noindent 
We linearize the system (\ref{sys1_1}) around a constant state $W_0$ defined by~:

\smallskip \centerline { $ \displaystyle   
W_{0}= ( \rho_0, u_0, v_0, s_0 )^t~,\, $ }

\smallskip  \noindent 
where $(...)^t$ is the transpose of a vector. The global state
 $\check{W}$ of the system is defined around the state $W_0$ thanks to the 
perturbation $\, W = ( \rho, u, v, s )^t \,$, as~:

\smallskip \centerline { $ \displaystyle  
\check{W} = W_0 + W \,$  . }

\smallskip  \noindent 
A first idea of our approach is to use the impulses~:
\[
\left \lbrace
  \begin{array}{l}
        \displaystyle{ \check{\rho} \check{u} = (\rho_0 + \rho) (u_0 + u) \equiv \rho_0 u_0 + \xi + \rho u} \\
        \displaystyle{ \check{\rho} \check{v} = (\rho_0 + \rho) (v_0 + v) \equiv \rho_0 v_0 + \zeta + \rho v}~,
  \end{array}
\right.
\]

\noindent and to linearize them considering the variables $\, \rho, u, v, s \,$ 
as first order infinitesimal quantities. We then introduce the {\bf{linearized impulses}}~:
\[
\left \lbrace
  \begin{array}{l}
  \xi = \rho_0 u + \rho u_0 \\
  \zeta = \rho_0 v + \rho v_0~.
  \end{array}
\right.
\]
\noindent We have the following classical hypothesis, see \cite{landau}~:\\
\vspace{-0.2cm}

\noindent {\textsc{{H}ypothesis 1}} \hspace{0.1cm} {\it Isentropy of the flow}.\\
\noindent  The linearization of the fourth equation of the system (\ref{sys1_1}) gives~:
\[
	\displaystyle{ \frac{\partial s}{\partial t} + u_0 \frac{\partial s}
{\partial x} + v_0 \frac{\partial s}{\partial y} \equiv \frac{ \rm{d} s}{ \rm{dt}} = 0 }~.
\]
\noindent If we consider the perturbation of entropy at the initial time to be null, 
that is to say $s(x,y,t=0) \equiv 0$, we deduce that $s(x,y,t) \equiv 0$ during the time 
evolution.

\smallskip  \noindent 
Then the system (\ref{sys1_1}) can be shared, first into a {\bf{stationary 
aerodynamic}} system~:
\[
  \left \lbrace
    \begin{array}{l}
      \mbox{div}~(\rho_{0} {\textbf{u}}_0) = 0 \\
      \rho_0 \textbf{u}_0 \bullet \nabla {\textbf{u}}_0 + \nabla p_0 = {\textbf{0}}~,
    \end{array}
  \right.
\]
\noindent and then into an {\bf{isentropic acoustic system}}~:

\moneq
\left \lbrace
  \begin{array}{l}
    \displaystyle{ \frac{\partial \rho}{\partial t}} + \frac{\partial \xi}
{\partial x} + \frac{\partial \zeta}{\partial y} = 0 \\
    \displaystyle{ \frac{\partial \xi}{\partial t} + \frac{\partial}{\partial x} 
\left( 2 u_0 \xi + \frac{c_0^2 - u_0^2}{c_0^2} p \right) + \frac{\partial}{\partial y} 
( u_0 \zeta + v_0 \xi - \rho u_0 v_0 ) = 0 } \\
    \displaystyle{ \frac{\partial \zeta}{\partial t} + \frac{\partial}{\partial x} 
( u_0 \zeta + v_0 \xi - \rho u_0 v_0 ) + \frac{\partial}{\partial y} \left( 2 v_0 \zeta 
+ \frac{c_0^2 - v_0^2}{c_0^2} p \right) = 0 }~,
  \end{array}
\right.
\label{sys1_3}
\end{equation}

\noindent with $\, p = c_0^2 \rho \,$ and $\, c_0 \,$ the speed of sound, deduced 
from the linearization of (\ref{sys1_2}). \\

\begin{proposition} [Advection of the acoustic vorticity] 
\end{proposition}

\noindent
If the external flow $W_{0}= ( \rho_0, u_0, v_0, s_0 )^t$ is stationary and uniform, 
the acoustic vorticity $\displaystyle{ \omega = \frac{\partial u}{\partial y} - 
\frac{\partial v}{\partial x} }$ is advected by the flow, {\it{i.e.}} 
$\displaystyle{ \frac{\rm{d \omega}}{\rm{dt}} = 0 }$.

\begin{proof}
This property is classical, see \cite{landau} for example. We give the 
proof for completeness. We have from the system (\ref{sys1_3})~:

\[
\begin{array}{lcl}
	\displaystyle{ \frac{\partial \xi}{\partial t} } & + & 
\displaystyle{ \frac{\partial}{\partial x} \left( 2 u_0 \xi + \frac{c_0^2 - u_0^2}{c_0^2} 
p \right) + \frac{\partial}{\partial y} ( u_0 \zeta + v_0 \xi - \rho u_0 v_0 ) } \\
	 & = & \displaystyle{ \frac{\partial \xi}{\partial t} + u_0 \frac{\partial}
{\partial x} \left( \xi - \rho u_0 \right) + v_0 \frac{\partial}{\partial y} ( \xi - \rho u_0) 
+ u_0 \frac{\partial \xi}{\partial x} + u_0 \frac{\partial \zeta}{\partial y} 
+ \frac{\partial p}{\partial x} } \\
	 & = & \displaystyle{ (\frac{\partial}{\partial t} + u_0 \frac{\partial}
{\partial x} + v_0 \frac{\partial}{\partial y}) (\xi - \rho u_0) + \frac{\partial p}{\partial x} = 0~,} \\
	\displaystyle{ \frac{\partial \zeta}{\partial t} } & + & \displaystyle{ 
\frac{\partial}{\partial x} ( u_0 \zeta + v_0 \xi - \rho u_0 v_0 ) + \frac{\partial}
{\partial y} \left( 2 v_0 \zeta + \frac{c_0^2 - v_0^2}{c_0^2} p \right) } \\
	 & = & \displaystyle{ \frac{\partial \zeta}{\partial t} + u_0 \frac{\partial}
{\partial x} \left( \zeta - \rho v_0 \right) + v_0 \frac{\partial}{\partial y} 
( \zeta - \rho v_0) + v_0 \frac{\partial \xi}{\partial x} + v_0 \frac{\partial \zeta}
{\partial y} + \frac{\partial p}{\partial y} } \\
	 & = & \displaystyle{ (\frac{\partial}{\partial t} + u_0 \frac{\partial}
{\partial x} + v_0 \frac{\partial}{\partial y}) (\zeta - \rho v_0) + \frac{\partial p}
{\partial y} = 0~.}
\end{array}
\]
\noindent We differentiate the first set of equations by $y$ and the second by $x$, 
we eliminate the pressure field and obtain~:

\moneq
	\displaystyle{ \frac{1}{\rho_0} \frac{\rm{d}}{\rm{dt}} \left( \frac{\partial}
{\partial y} ( \xi - \rho u_0 ) - \frac{\partial}{\partial x} ( \zeta - \rho v_0 ) \right) 
= \frac{\rm{d}}{\rm{dt}} ( \frac{\partial u}{\partial y} - \frac{\partial v}{\partial x} ) 
= \frac{\rm{d \omega}}{\rm{dt}} = 0 }~.
\label{sys1_4}
\end{equation}
\end{proof}

\noindent {\textsc{{H}ypothesis 2}} \hspace{0.1cm}  { \it Irrotationality 
of the acoustic vorticity.} 

\noindent If we consider the acoustic perturbation at the initial time to be 
{\bf{irrotational}}, {\it{i.e.}} ${\bf{rot\,u}} \, (x,y,t=0) \equiv {\bf{\omega}} \, 
(x,y,t=0) = {\bf{0}}$, we then deduce with equation (\ref{sys1_4}) that 
${\bf{rot\,u}} \, (x,y,t) = {\bf{0}}$ during the time evolution.

\bigskip \bigskip  \noindent {\bf \large 3) \quad Lorentz Transform} 

\noindent 
We consider the two-dimensional equations of advective acoustics when the velocity 
of the fluid is parallel to a particular direction; we suppose specifically~:

\moneq
	{\bf{u}} = u_0~{\bf{e_x}}~.
\label{sys2_1}
\end{equation}

\noindent We search a space-time transform $(x,t) \longmapsto (x',t')$ so that in 
the new space-time $(x',t')$, the pressure field is the solution of the wave equation. 
We find that this space-time transform is a {\bf{Lorentz transform}}. With it, we 
derive a new set of equations and prove that the corresponding system can be reduced 
to the classical case of non advective acoustics.

\bigskip    \noindent {\bf \large 3-1  \quad Change of space-time} 

\noindent 
Considering a flow of velocity given by equation (\ref{sys2_1}), the system (\ref{sys1_3}) is written as~:

\moneq
\left \lbrace
  \begin{array}{l}
    \displaystyle{ \frac{\partial p}{\partial t} + c_0^2 \frac{\partial \xi}{\partial x} 
+ c_0^2 \frac{\partial \zeta}{\partial y} = 0 } \\
    \displaystyle{ \frac{\partial \xi}{\partial t} + \frac{\partial}{\partial x} 
\left( 2 u_0 \xi + \frac{(c_0^2 - u_0^2)}{c_0^2} p \right) + \frac{\partial}{\partial y} 
( u_0 \zeta ) = 0 } \\
    \displaystyle{ \frac{\partial \zeta}{\partial t} + \frac{\partial}{\partial x} 
( u_0 \zeta ) + \frac{\partial p}{\partial y } = 0 }~,
  \end{array}
\right.
\label{sys2_2} 
\end{equation}

\noindent which is a pleasant conservative form. We easily deduce that the pressure 
field $p \, (x,y,t)$ is solution in the (initial) space-time $(x,y,t)$ of a wave
equation~:

\moneq
        \displaystyle{ \frac{\partial^2 p}{\partial t^2} + 2 u_0 \frac{\partial^2 p}
{\partial x \partial t} + u_0^2 \frac{\partial^2 p}{\partial x^2} - c_0^2 \Delta p = 0 }~,
\label{sys2_3}
\end{equation}
\noindent where $\displaystyle{ \Delta = \frac{\partial^2}{\partial x^2} +
  \frac{\partial^2}
{\partial y^2}}$ is the laplacian in two dimension space. \\

\begin{proposition}[Lorentz transform and equation of pressure]
\label{prop2_1}
\end{proposition}

\nototo
We suppose that the advective velocity satisfies equation (\ref{sys2_1}). We define the 
Mach number as $ {M_0 = \frac{u_0}{c_0}}$ and the Lorentz space-time transform
as~:

\moneq
\left \lbrace
  \begin{array}{l}
        x' = \displaystyle{ \frac{1}{\sqrt{1-M_0^2}} } \, x \\
        y' = y \\
        t' = t + \displaystyle{ \frac{M_0}{c_0 (1-M_0^2)} } \, x~.
  \end{array}
\right.
\label{sys2_4}
\end{equation}
\noindent In this new space-time, the pressure field is considered as a function of 
the new set of space-time coordinates $(x',y',t')$, {\it i.e.}~:

\moneq
	p'(x',y',t') \equiv p \, (x,y,t)~,
\label{sys2_5}
\end{equation}
\noindent and is the solution of the wave equation with a {\bf{modified celerity}}~:

\moneq
	\frac{\partial^2 p'}{\partial {t'}^2} - c_0^2 (1-M_0^2) \left( 
\frac{\partial^2 p'}{\partial {x'}^2} + \frac{\partial^2 p'}{\partial {y'}^2} \right) = 0~.
\label{sys2_6}
\end{equation}
\noindent reduced from the ``pure'' sound celerity by a similarity factor $\sqrt{1-M_0^2}$~.

\begin{proof}
We first explain the way we derive the Lorentz transform (\ref{sys2_4}) to remove 
the advective contribution $2 u_0 \frac{\partial^2 p}{\partial x \partial t}$ in 
equation (\ref{sys2_3}). In the new space-time $(x',y',t')$, we want the pressure 
field to be solution of the wave equation. We search the new space-time coordinates 
$(x',y',t')$ as~:

\moneq
   \left \lbrace
   \begin{array}{l}
        x' = \alpha x \\
        y' = y \\
        t' = t + \beta x~.
    \end{array}
\right.
\label{change_espace_temps}
\end{equation}

\noindent The transformed equation (\ref{sys2_3}) takes the form~:

\[
   \begin{array}{lcl}
        \displaystyle{ \frac{\partial^2 p}{\partial t^2} } & + & \displaystyle{ 2 u_0 
\frac{\partial^2 p}{\partial x \partial t} + u_0^2 \frac{\partial^2 p}{\partial x^2} 
- c_0^2 \Delta p} \\
         & = & \displaystyle{ \left[ \left( \frac{\partial}{\partial t} + u_0 
\frac{\partial}{\partial x} \right)^2 - c_0^2 \Delta \right] p\,(x,y,t) } \\
         & = & \displaystyle{ \left[ \left( (1+u_0 \beta)^2 -c_0^2 \beta^2 \right) 
\frac{\partial^2}{\partial {t^{'}}^2} + 2 \alpha \left( u_0 ( 1 + u_0 \beta ) - 
\beta c_0^2 \right) \frac{\partial^2}{\partial t' \partial x'} \right. } \\
         &  & \displaystyle{ \left. - \alpha^2 (c_0^2 - u_0^2) 
\frac{\partial^2}{\partial {x^{'}}^2} - c_0^2 \frac{\partial^2}{\partial {y^{'}}^2} \right] p'(x',y',t') }~.
   \end{array}
\]

\noindent Then we obtain~:

\moneq
   \begin{array}{l}
	\displaystyle{ \left[ \left( (1+u_0 \beta)^2 -c_0^2 \beta^2 \right)
 \frac{\partial^2}{\partial {t^{'}}^2} \right.}  +  \displaystyle{ \left. 2 \alpha 
\left( u_0 ( 1 + u_0 \beta ) - \beta c_0^2 \right) \frac{\partial^2}{\partial t' \partial
 x'} \right. } \\
\qquad \quad  -  \displaystyle{ \left. \alpha^2 (c_0^2 - u_0^2) 
\frac{\partial^2}{\partial {x^{'}}^2} - c_0^2 \frac{\partial^2}{\partial 
{y^{'}}^2} \right] p'(x',y',t') = 0 }~.
    \end{array}
\label{neweqpression}
\end{equation}

\noindent The conditions upon $\alpha$ and $\beta$ to find the wave equation are 
clear from equation (\ref{neweqpression}); on the first hand no further crossed
 partial derivation between space and time, that is~:

\moneq 
        \displaystyle{ 2 \alpha \left( u_0 ( 1 + u_0 \beta ) - \beta c_0^2 \right) = 0 }~,
\label{condition1}
\end{equation}

\noindent and on the other hand equality of the coefficients of double derivations in
 space to have a laplacian operator invariant by rotation~:

\moneq
        \displaystyle{ \alpha^2 (c_0^2 - u_0^2) = c_0^2 }~.
\label{condition2}
\end{equation}

\noindent The unique solution $(\alpha,\beta)$ of the previous 2\,x\,2 linear 
system (\ref{condition1})-(\ref{condition2}) is~:

\[
\left \lbrace
   \begin{array}{l}
        \displaystyle{ \alpha = \frac{c_0}{\sqrt{c_0^2 - u_0^2}} } \\
        \displaystyle{ \beta  = \frac{u_0}{c_0^2 - u_0^2} }~,
    \end{array}
\right.
\]

\noindent and with this set of coefficients, the space-time transform 
(\ref{change_espace_temps}) is exactly equal to the system (\ref{sys2_4}). Moreover, 
we remark that the coefficient of $\displaystyle{ \frac{\partial^2}{\partial t^2} }$ 
in equation (\ref{neweqpression}) is now equal to~:
\[ \displaystyle{ \left( (1+u_0 \beta)^2 -c_0^2 \beta^2 \right) = \frac{c_0^2}{c_0^2 
- u_0^2} = \frac{1}{1 - M_0^2} }~, \]
and, in our transformed space-time $(x',y',t')$, the pressure field defined by the 
condition (\ref{sys2_5}) is the solution of the wave equation (\ref{sys2_6})~.
\end{proof}

\bigskip    \noindent {\bf \large 3-2  \quad Change of unknown functions} 

\noindent 
Let us apply the Lorentz transform (\ref{sys2_4}) to the acoustic system 
(\ref{sys2_2}). We have the following proposition~:

\begin{proposition}[New unknown functions for advective acoustic]
\label{prop2_2}
\end{proposition}

\nototo
We assume that the hypothesis 2 of irrotationality of the acoustic vorticity is 
satisfied and that the advective velocity field is defined by equation (\ref{sys2_1}). 
After applying the Lorentz transform (\ref{sys2_4}) and the following change of 
pressure and impulse functions~: 

\moneq
\left \lbrace
  \begin{array}{l}
        \displaystyle { \widetilde{p} = p' + \frac{u_0}{(1-M_0^2)} \, \xi' }  \\
        \displaystyle { \widetilde{\xi} = \frac{1}{\sqrt{1-M_0^2}} \, \xi' } \\
        \displaystyle { \widetilde{\zeta} = \zeta' }~,
  \end{array}
\right.
\label{sys2_7}
\end{equation}

\noindent the advective acoustic system (\ref{sys2_2}) can be written as~:

\moneq
\left \lbrace
  \begin{array}{l}
        \displaystyle{ \frac{\partial \widetilde{p}}{\partial t'} + c_0^2 
\frac{\partial \widetilde{\xi}}{\partial x'} + c_0^2 \frac{\partial \widetilde{\zeta}}
{\partial y'} = 0 } \\
        \displaystyle{ \frac{\partial \widetilde{\xi}}{\partial t'} + (1-M_0^2) 
\frac{\partial \widetilde{p}}{\partial x'} = 0 } \\
        \displaystyle{ \frac{\partial \widetilde{\zeta}}{\partial t'} + (1-M_0^2) 
\frac{\partial \widetilde{p}}{\partial y'} = 0 }~.
  \end{array}
\right.
\label{sys2_8}
\end{equation}

\begin{proof}
We first use the hypothesis of irrotationality of the acoustic vorticity in 
the third equation of the system (\ref{sys2_2}) and obtain~: 
\[
	\displaystyle{ \frac{\partial \zeta}{\partial x} = \frac{\partial 
(\xi - \rho u_0)}{\partial y} = \frac{\partial \xi}{\partial y} - \frac{u_0}{c_0^2} 
\frac{\partial p}{\partial y} }~.
\]
\noindent Secondly we introduce the Lorentz transform (\ref{sys2_4}) into the 
system (\ref{sys2_2}). We have the following transform of partial derivations~:
\[
    \left \lbrace
      \begin{array}{l}
	\displaystyle{ \frac{\partial}{\partial x} = \frac{1}{\sqrt{1-M_0^2}} 
\frac{\partial}{\partial x'} + \frac{u_0}{c_0^2 (1-M_0^2) } \frac{\partial}{\partial t'} } \\
	\displaystyle{ \frac{\partial}{\partial y} = \frac{\partial}{\partial y'} } \\
	\displaystyle{ \frac{\partial}{\partial t} = \frac{\partial}{\partial t'} }~.
      \end{array}
    \right.
\]
\noindent We then substract the first equation of the system (\ref{sys2_2})
 multiplied by $\, \frac{u_0}{c_0^2} \,$ from the second one and, using the following
 notations~:

\moneq
 \left \lbrace
    \begin{array}{l}
	p'(x',y',t') \equiv p\,(x,y,t) \\
	\xi'(x',y',t') \equiv \xi (x,y,t) \\
	\zeta'(x',y',t') \equiv \zeta (x,y,t)~,
    \end{array}
 \right.
\label{fonctionsprime}
\end{equation}
\noindent we find~:
\[
\left \lbrace
  \begin{array}{l}
        \displaystyle{ \frac{\partial p'}{\partial t'} + \frac{c_0^2}{\sqrt{1-M_0^2}} 
\frac{\partial \xi'}{\partial x'} + \frac{M_0 c_0}{(1-M_0^2)} \frac{\partial \xi'}
{\partial t'} + c_0^2 \frac{\partial \zeta'}{\partial y'} = 0 } \\
        \displaystyle{ \frac{\partial \xi'}{\partial t'} + \frac{u_0}{\sqrt{1-M_0^2}} 
\frac{\partial \xi'}{\partial x'} + \frac{M_0^2}{1-M_0^2} \frac{\partial \xi'}{\partial t'} 
+ \frac{c_0^2 - u_0^2}{c_0^2 \sqrt{1-M_0^2}} \frac{\partial p'}{\partial x'} = 0 } \\
        \displaystyle{\frac{ \partial \zeta'}{\partial t'} + u_0 \frac{\partial \xi'}
{\partial y'} + \frac{( c_0^2 - u_0^2 )}{c_0^2} \frac{\partial p'}{\partial y'} = 0 }~.
  \end{array}
\right.
\]
\noindent We gather the terms associated with the same operator of derivation~:
\[
\left \lbrace
  \begin{array}{l}
        \displaystyle{ \frac{\partial}{\partial t'} \left[ p' + \frac{M_0 c_0}{(1-M_0^2)} 
\xi' \right] + c_0^2 \frac{\partial}{\partial x'} (\frac{\xi'}{\sqrt{1-M_0^2}}) + c_0^2 
\frac{\partial \zeta'}{\partial y'} = 0 } \\
        \displaystyle{ \frac{\partial}{\partial t'} (\frac{\xi'}{\sqrt{1-M_0^2}}) + 
(1-M_0^2) \frac{\partial}{\partial x'} \left[ p' + \frac{M_0 c_0}{(1-M_0^2)} \xi' 
\right] = 0 } \\
        \displaystyle{ \frac{\partial \zeta'}{\partial t'} + (1-M_0^2) \frac{\partial}
{\partial y'} \left[ p' + \frac{M_0 c_0}{(1-M_0^2)} \xi' \right] = 0 }~,
  \end{array}
\right.
\]
\noindent and we substitute into the previous system the new unknown functions 
$(\widetilde{p} \, , \widetilde{\xi} \, , \widetilde{\zeta})$ introduced in the system 
(\ref{sys2_7}). Then the system of equations (\ref{sys2_8}) is satisfied.
\end{proof}

\smallskip \noindent {\bf   Remark 4. } 

\noindent 
The major consequence of propositions \ref{prop2_1} and \ref{prop2_2} 
is the following (operational\,!) remark. The resolution of the advective 
acoustic system is absolutly {\bf{identical}} to the one obtained {\bf{without}} 
advective flow, but with a propagation celerity {\bf{scaled}} by a factor
$\sqrt{1-M_0^2}$.

\bigskip \bigskip  \noindent {\bf \large 4) \quad Lorentz transform for multi-dimensional flows}

\noindent 
In the previous section, we dealt with the case of a velocity field described 
by the equation (\ref{sys2_1}). In \cite{gottlieb}, Abarbanel {\it{et}} al consider 
a multi-dimensional flow as a one-dimensional flow, after a correct rotation of the 
studied medium by an angle $\theta = \tan^{-1} (\frac{v_0}{u_0})$ and considering
 the new velocity to be $u_{new} = \sqrt{u_0^2 + v_0^2}$. We observe that in order 
to study numerically the influence of the flow for acoustic propagation near objects, 
such an idea imposes a remeshing of the geometry for each change of the advective flow. 
In our opinion, this process is not compatible with the use of finite differences and 
with operational industrial constraints. We propose in this section to generalize 
the Lorentz space-time transform to a multi-dimensional flow and to extend our approach 
with the help of {\bf{space affinities}} to a multi-dimensional flow under the same
 hypotheses as before. In the next two paragraphs, we present the generalization of 
the Lorentz transform respectively to the two and three-dimensional cases. Only the 
two-dimensional case is proven in the present document. The proof of the 
three-dimensional case can be found in \cite{fred}.

\bigskip    \noindent {\bf \large 4-1  \quad The two-dimensional case}  

\noindent 
We now consider a subsonic uniform flow described by a velocity vector~:

\moneq
	{\bf{u}} = (u_0,v_0)~.
\label{sys2_9}
\end{equation}

\noindent With such an external flow, the linearized isentropic Euler equations 
for advective acoustics are~:

\moneq
  \left \lbrace
    \begin{array}{l}
      \displaystyle{ \frac{\partial p}{\partial t}} + c_0^2 \frac{\partial \xi}{\partial x} 
+ c_0^2 \frac{\partial \zeta}{\partial y} = 0 \\
      \displaystyle{ \frac{\partial \xi}{\partial t} + \frac{\partial}{\partial x} 
\left( 2 u_0 \xi + \frac{(c_0^2 - u_0^2)}{c_0^2} p \right) + \frac{\partial}{\partial y}
 ( u_0 \zeta + v_0 \xi - \frac{u_0 v_0}{c_0^2} p ) = 0 } \\
      \displaystyle{ \frac{\partial \zeta}{\partial t} + \frac{\partial}{\partial x} 
( u_0 \zeta + v_0 \xi - \frac{u_0 v_0}{c_0^2} p ) + \frac{\partial}{\partial y} 
\left( 2 v_0 \zeta + \frac{(c_0^2 - v_0^2)}{c_0^2} p \right) = 0 }~.
    \end{array}
  \right.
\label{sys2_10}
\end{equation}

\begin{proposition}[Simplification of the acoustic system]
\end{proposition}

\nototo
We generalize the Lorentz space-time transform when the external 
flow verifies (\ref{sys2_9}). We introduce a new set of space-time coordinates as~:

\moneq
  \left \lbrace
    \begin{array}{l}
        x' = \displaystyle{ \frac{1}{\sqrt{1-\frac{u_0^2}{c_0^2}}} \, x } \\
        y' = \displaystyle{ \frac{1}{\sqrt{1-\frac{v_0^2}{c_0^2}}} \, y } \\
        t' = t + \displaystyle{ \frac{u_0}{c_0^2 (1-M_0^2)} } \, x + 
\displaystyle{ \frac{v_0}{c_0^2 (1-M_0^2)} } \, y~,
    \end{array}
  \right.
\label{sys2_11}
\end{equation}

\noindent the Mach number as
 $M_0 = \frac{\sqrt{u_0^2 + v_0^2}}{c_0}$, a coupling coefficient $\alpha$ 
between the two cartesian coordinates as~:

\moneq
        \alpha = \frac{u_0 v_0}{c_0^2 \sqrt{1-\frac{u_0^2}{c_0^2}} \sqrt{1-\frac{v_0^2}{c_0^2}}}~,
\label{sys2_13}
\end{equation}

\noindent and the new unknown functions $\widetilde{p}$, $\widetilde{\xi}$ et $\widetilde{\zeta}$ 
defined by~:

\moneq
\left \lbrace
  \begin{array}{l}
        \displaystyle{ \widetilde{p} = p' + \frac{1}{1-M_0^2} \left( u_0 \xi' 
+ v_0 \zeta' \right) } \\ 
        \displaystyle{ \widetilde{\xi} = \sqrt{ 1-\frac{u_0^2}{c_0^2} } \left( 
(1 - \frac{v_0^2}{c_0^2}) \frac{\xi'}{1-M_0^2} + \frac{u_0 v_0}{c_0^2} 
\frac{\zeta'}{1-M_0^2} \right) } \\
        \displaystyle{ \widetilde{\zeta} = \sqrt{ 1-\frac{v_0^2}{c_0^2} } 
\left( \frac{u_0 v_0}{c_0^2} \frac{\xi'}{1-M_0^2} + (1 - \frac{u_0^2}{c_0^2}) 
\frac{\zeta'}{1-M_0^2} \right) }~.
  \end{array}
\right.
\label{sys2_14}
\end{equation}

\noindent Under the hypothesis 2 of irrotationality of the acoustic vorticity, 
the new form of the acoustic system (\ref{sys2_10}) is the following~:

\moneq
  \left \lbrace
    \begin{array}{l}
        \displaystyle{ \frac{\partial \widetilde{p}}{\partial t'} + c_0^2 
\frac{\partial}{\partial x'} ( \widetilde{\xi} - \alpha \widetilde{\zeta} ) + c_0^2 
\frac{\partial}{\partial y'} ( \widetilde{\zeta} - \alpha \widetilde{\xi} ) = 0 } \\
        \displaystyle{ \frac{\partial \widetilde{\xi}}{\partial t'} + (1-M_0^2) 
\frac{\partial \widetilde{p}}{\partial x'} = 0 } \\
        \displaystyle{ \frac{\partial \widetilde{\zeta}}{\partial t'} + (1-M_0^2) 
\frac{\partial \widetilde{p}}{\partial y'} = 0 }~.
    \end{array}
  \right.
\label{sys2_12}
\end{equation}

\begin{proof} 
First, by substracting with a correct coefficient the first equation of the system 
(\ref{sys2_10}) from the two others, we find~:

\[
  \left \lbrace
    \begin{array}{l}
      \displaystyle{ \frac{\partial p}{\partial t}} + c_0^2 
\frac{\partial \xi}{\partial x} + c_0^2 \frac{\partial \zeta}{\partial y} = 0 \\
      \displaystyle{ \frac{\partial}{\partial t} ( \xi - \frac{u_0}{c_0^2} p ) 
+ \frac{\partial}{\partial x} \left( u_0 \xi + \frac{(c_0^2 - u_0^2)}{c_0^2} p 
\right) + \frac{\partial}{\partial y} ( v_0 \xi - \frac{u_0 v_0}{c_0^2} p ) = 0 } \\
      \displaystyle{ \frac{\partial}{\partial t} ( \zeta - \frac{v_0}{c_0^2} p ) 
+ \frac{\partial}{\partial x} ( u_0 \zeta - \frac{u_0 v_0}{c_0^2} p ) + 
\frac{\partial}{\partial y} \left( v_0 \zeta + \frac{(c_0^2 - v_0^2)}{c_0^2} p \right) = 0 }~.
    \end{array}
  \right.
\]

\noindent Using the hypothesis of irrotationality of the acoustic vorticity, we have~:

\[
   \begin{array}{lcl}
        \displaystyle{ \frac{\partial}{\partial y} (v_0 \xi - 
\frac{u_0 v_0}{c_0^2} p) } & = & \displaystyle{ \frac{\partial}{\partial y} 
\left[ v_0 (\rho_0 u + \rho u_0) - \frac{u_0 v_0}{c_0^2} p \right] = 
\frac{\partial}{\partial y} (\rho_0 v_0 u) } \\
         & = & \displaystyle{ \frac{\partial}{\partial x} (\rho_0 v_0 v) = 
\frac{\partial}{\partial x} v_0 ( \zeta - \rho v_0 ) } \\
         & = & \displaystyle{ \frac{\partial}{\partial x} ( v_0 \zeta - 
\frac{v_0^2}{c_0^2} p ) }~.
   \end{array}
\]

\noindent The previous calculation gives the new system~:

\[
  \left \lbrace
    \begin{array}{l}
      \displaystyle{ \frac{\partial p}{\partial t}} + c_0^2 
\frac{\partial \xi}{\partial x} + c_0^2 \frac{\partial \zeta}{\partial y} = 0 \\
      \displaystyle{ \frac{\partial}{\partial t} ( \xi - \frac{u_0}{c_0^2} p ) 
+ \frac{\partial}{\partial x} \left[ u_0 \xi + v_0 \zeta + 
(1-M_0^2) p \, \right]  = 0 } \\
      \displaystyle{ \frac{\partial}{\partial t} ( \zeta - \frac{v_0}{c_0^2} p )
 + \frac{\partial}{\partial y} \left[ u_0 \xi + v_0 \zeta + (1-M_0^2) p \, 
\right] = 0 }~.
    \end{array}
  \right.
\]

\noindent In the new space-time defined by the change of space-time 
(\ref{sys2_11}), we have~:
\[
    \left \lbrace
      \begin{array}{l}
	\displaystyle{ \frac{\partial}{\partial x} = 
\frac{1}{\sqrt{1-\frac{u_0^2}{c_0^2}}} \frac{\partial}{\partial x'} 
+ \frac{u_0}{c_0^2 (1-M_0^2) } \frac{\partial}{\partial t'} } \\
	\displaystyle{ \frac{\partial}{\partial y} = 
\frac{1}{\sqrt{1-\frac{v_0^2}{c_0^2}}} \frac{\partial}{\partial y'} +
 \frac{v_0}{c_0^2 (1-M_0^2) } \frac{\partial}{\partial t'} } \\
	\displaystyle{ \frac{\partial}{\partial t} = 
\frac{\partial}{\partial t'} }~.
      \end{array}
    \right.
\]

\noindent With the notation introduced in the system (\ref{fonctionsprime}), 
the transformed equations take the algrebraic form~:

\[
  \left \lbrace
    \begin{array}{l}
      \displaystyle{ \frac{\partial}{\partial t'}} \left( p' + \frac{u_0 \xi' 
+ v_0 \zeta'}{1-M_0^2} \right) + c_0^2 \frac{\partial}{\partial x'} ( 
\frac{\xi'}{\sqrt{1-\frac{u_0^2}{c_0^2}}} ) + c_0^2 \frac{\partial}{\partial y'} 
( \frac{\zeta'}{\sqrt{1-\frac{v_0^2}{c_0^2}}} ) = 0 \\
      \displaystyle{ \frac{\partial}{\partial t'} \left[ \frac{\sqrt{
 1-\frac{u_0^2}{c_0^2}}}{1-M_0^2} \left( (1 - \frac{v_0^2}{c_0^2}) \xi' + 
\frac{u_0 v_0}{c_0^2} \zeta' \right) \right] + (1-M_0^2) \frac{\partial}{\partial x'} 
\left( p' + \frac{u_0 \xi' + v_0 \zeta'}{1-M_0^2} \right) = 0 } \\
      \displaystyle{ \frac{\partial}{\partial t'} \left[ \frac{\sqrt{ 
1-\frac{v_0^2}{c_0^2}}}{1-M_0^2} \left( \frac{u_0 v_0}{c_0^2} \xi'  + 
(1 - \frac{u_0^2}{c_0^2}) \zeta' \right) \right] + (1-M_0^2) \frac{\partial}
{\partial y'} \left( p' + \frac{u_0 \xi' + v_0 \zeta'}{1-M_0^2} \right) = 0 }~.
    \end{array}
  \right.
\]

\noindent By the change of unknown functions (\ref{sys2_14}), the previous system
becomes~:

\moneq
\left \lbrace
  \begin{array}{l}
    \displaystyle{ \frac{\partial \widetilde{p}}{\partial t'}} + c_0^2 
\frac{\partial}{\partial x'} ( \frac{\xi'}{\sqrt{1-\frac{u_0^2}{c_0^2}}} )
 + c_0^2 \frac{\partial}{\partial y'} ( \frac{\zeta'}{\sqrt{1-\frac{v_0^2}{c_0^2}}} ) = 0 \\
    \displaystyle{ \frac{\partial \widetilde{\xi}}{\partial t'} + (1-M_0^2) 
\frac{\partial \widetilde{p}}{\partial x'} = 0 } \\
    \displaystyle{ \frac{\partial \widetilde{\zeta}}{\partial t'} + (1-M_0^2) 
\frac{\partial \widetilde{p}}{\partial y'} = 0 }~.
  \end{array}
\right.
\label{sysxiprime}
\end{equation}

\noindent We focus here on the fact that the pair $(\xi', \zeta')$ is present 
in the first equation of (\ref{sysxiprime}) whereas the new unknown functions 
are $\widetilde{\xi}$ and $\widetilde{\zeta}$. Nevertheless with the last two 
equations of the system (\ref{sys2_14}), we have the following calculation~:

\[
\left \lbrace
    \begin{array}{l}
	\displaystyle{ \xi' = \sqrt{1 - \frac{u_0^2}{c_0^2}}~( \widetilde{\xi} 
- \alpha \widetilde{\zeta} ) } \\
	\displaystyle{ \zeta' = \sqrt{1 - \frac{v_0^2}{c_0^2}}~( 
\widetilde{\zeta} - \alpha \widetilde{\xi} ) }~,
    \end{array}
\right.
\]

\noindent where $\alpha$ is defined by equation (\ref{sys2_13}). We then find 
the final form (\ref{sys2_12}) of the system of advective acoustics in the new 
space-time $(x',y',t')$.
\end{proof}

\bigskip    \noindent {\bf \large 4-2 \quad The three-dimensional case} 

\noindent 
The generalization to {\bf{three dimension space}} can be done without 
any major difficulty. We have the following proposition proven in \cite{fred}~:

\begin{proposition} 
\end{proposition}

\nototo
We assume that the velocity of the external flow is given by~:

\moneq
	\displaystyle{ {\bf{u_0}} = (u_0, v_0, w_0) }~,
\end{equation}

\noindent and that the Hypothesis 2 of irrotationality of the 
acoustic vorticity is verified. We define the Mach number as $M_0 = 
\frac{\sqrt{u_0^2 + v_0^2 + w_0^2}}{c_0}$ and a change of space-time as~:

\moneq
\left \lbrace
  \begin{array}{l}
        \displaystyle{ x' = \frac{1}{\sqrt{1-\frac{u_0^2}{c_0^2}}} \, x } \\
        \displaystyle{ y' = \frac{1}{\sqrt{1-\frac{v_0^2}{c_0^2}}} \, y } \\
        \displaystyle{ z' = \frac{1}{\sqrt{1-\frac{w_0^2}{c_0^2}}} \, z } \\
        \displaystyle{ t' = t+ \frac{u_0}{c_0^2 (1-M_0^2)} \, x + 
\frac{v_0}{c_0^2 (1-M_0^2)} \, y + \frac{w_0}{c_0^2 (1-M_0^2)} \, z }~.
  \end{array}
\right.
\end{equation}

\noindent We introduce {\bf{three}} coupling coefficients as~: 

\[
  \left \lbrace
     \begin{array}{lcl}
        \displaystyle{ \alpha } & = & \displaystyle{ \frac{u_0 v_0}{c_0^2 
\sqrt{1-\frac{u_0^2}{c_0^2}} \sqrt{1-\frac{v_0^2}{c_0^2}}} } \\
        \displaystyle{ \beta } & = & \displaystyle{ \frac{u_0 w_0}{c_0^2 
\sqrt{1-\frac{u_0^2}{c_0^2}} \sqrt{1-\frac{w_0^2}{c_0^2}}} } \\
        \displaystyle{ \gamma } & = & \displaystyle{ \frac{v_0 w_0}{c_0^2 
\sqrt{1-\frac{v_0^2}{c_0^2}} \sqrt{1-\frac{w_0^2}{c_0^2}}} }~,
     \end{array}
  \right.
\]

\noindent and a change of unknown functions as~:

\[
\left \lbrace
  \begin{array}{l}
        \displaystyle{ \widetilde{p} = p' + \frac{1}{1-M_0^2} 
\left( u_0 \xi' + v_0 \zeta' + w_0 \chi' \right) } \\ 
        \displaystyle{ \widetilde{\xi} = \sqrt{ 1-\frac{u_0^2}{c_0^2} } 
\left( (1 - \frac{v_0^2+w_0^2}{c_0^2}) \frac{\xi'}{1-M_0^2} + \frac{u_0 v_0}{c_0^2} 
\frac{\zeta'}{1-M_0^2} + \frac{u_0 w_0}{c_0^2} \frac{\chi'}{1-M_0^2} \right) } \\
        \displaystyle{ \widetilde{\zeta} = \sqrt{ 1-\frac{v_0^2}{c_0^2} } 
\left( \frac{u_0 v_0}{c_0^2} \frac{\xi'}{1-M_0^2} + (1 - \frac{u_0^2+w_0^2}{c_0^2}) 
\frac{\zeta'}{1-M_0^2} + \frac{v_0 w_0}{c_0^2} \frac{\chi'}{1-M_0^2} \right) } \\
        \displaystyle{ \widetilde{\chi} = \sqrt{ 1-\frac{w_0^2}{c_0^2} } 
\left( \frac{u_0 w_0}{c_0^2} \frac{\xi'}{1-M_0^2} + \frac{v_0 w_0}{c_0^2} 
\frac{\zeta'}{1-M_0^2} + (1 - \frac{u_0^2+v_0^2}{c_0^2}) \frac{\chi'}{1-M_0^2} \right) }~.
  \end{array}
\right.
\]

\noindent The acoustic system takes the new form~:

\[
\left \lbrace
  \begin{array}{l}
      \begin{array}{lcl}
        \displaystyle{ \frac{\partial \widetilde{p}}{\partial t'} } & + & 
\displaystyle{ c_0^2 \frac{\partial}{\partial x'} ( \widetilde{\xi} - \alpha 
\widetilde{\zeta} - \beta \widetilde{\chi} ) } \\
         & + & \displaystyle{ c_0^2 \frac{\partial}{\partial y'} 
( \widetilde{\zeta} - \alpha \widetilde{\xi} - \gamma \widetilde{\chi} ) + 
c_0^2 \frac{\partial}{\partial z'} ( \widetilde{\chi} - \beta \widetilde{\xi} - 
\gamma \widetilde{\zeta} ) = 0 }
      \end{array} \\
        \displaystyle{ \frac{\partial \widetilde{\xi}}{\partial t'} + (1-M_0^2) 
\frac{\partial \widetilde{p}}{\partial x'} = 0 } \\
        \displaystyle{ \frac{\partial \widetilde{\zeta}}{\partial t'} + (1-M_0^2) 
\frac{\partial \widetilde{p}}{\partial y'} = 0 } \\
        \displaystyle{ \frac{\partial \widetilde{\chi}}{\partial t'} + (1-M_0^2) 
\frac{\partial \widetilde{p}}{\partial z'} = 0 }~.
  \end{array}
\right.
\]

\bigskip \bigskip  \noindent {\bf \large 5) \quad Acoustic Absorbing Layers} 

\noindent  
Physical wave phenomena modelling takes often place in the infinite two or 
three dimensional space. Due to finite computing resources, the numerical 
simulations of such phenomena must be truncated to confined domains, then numerical 
artifical boundaries must be considered. Generally, numerical reflections of 
outgoing waves from the boundaries of the numerical domain reenter the computational
domain and falsify the results. Various methods have been proposed to reduce the influence
of reentering waves in the 
computational domain.

\smallskip \noindent 
For many years, the numerical physicists, Israeli-Orsag \cite{israeli}, 
have developed the idea of layers of absorbing materials. Then the mathematical 
study of non reflecting boundary conditions has been developed after the pionnering 
work of Engquist-Majda \cite{engquist}. The discrete studies of such absorbing 
conditions have been realized for scalar waves by Bayliss-Turkel \cite{bayliss}, 
for electromagnetic waves by Joly-Mercier \cite{joly}, Taflove \cite{taflove} 
and for sismic waves by Halpern-Trefethen \cite{halpern} among others.

\smallskip \noindent 
The current perfectly matched layers approach has been introduced by 
B\'erenger \cite{berenger} in the context of computational electromagnetics; 
a mathematical interpretation of this model has been made by Collino \cite{collino}. 
In \cite{hu}, Hu proposes an adaptation of Berenger's model for advective acoustics. 
Nevertheless, Abarbanel {\it{et}} al \cite{gottlieb} and Rahmouni \cite{adib2} and 
\cite{adib1} have demonstrated that this model is mathematically ill-posed, {\it{i.e.}} 
that, if truncated to the first order terms, there exists a perturbation as small as
 we wish that can make the model unstable. These authors also propose well-posed models.

\smallskip \noindent 
Our approach uses a model of the type introduced by Hu. We focus in our study 
on the practical disadvantages to deal with a mathematical model whose principal
 symbol corresponds to a ill-posed problem.

\bigskip  \newpage   \noindent {\bf \large 5-1  \quad Acoustic absorbing layers without external flow}  

\noindent 
We propose in the following a precise description of acoustic absorbing layers, 
and we follow the ideas developed by Collino \cite{collino}. We consider a 
semi-infinite medium in the $x$-direction defined by $\Omega = \Omega^{+} \cup 
\Omega^{-}$, where~:

\moneq
\left \lbrace
	\begin{array}{ll}
		\Omega^{-} = \{ (x,y),\mbox{~}y<0 \} \\
		\Omega^{+}= \{ (x,y),\mbox{~}0 \le y \le \delta \}~,
	\end{array}
\right.
\label{def_omega_plus}
\end{equation}

\noindent $\Omega^{+}$ representing the absorbing layers domain.\\

\begin{center}
   \begin{picture}(60,52)
	\put(25,0){\circle*{1.5}}
	\qbezier[100](21,0)(25,6)(29,0)
	\qbezier[200](17,0)(25,12)(33,0)
	\qbezier[300](12,0)(25,18)(37,0)
	\put(15,15){\large{$\Omega^{-}$}}
	\put(32,18){Physical}
	\put(33,13){domain}
	\put(0,24){\line(1,0){50}}
	\put(52,23){$y=0$}
	\put(15,30){\large{$\Omega^{+}$}}
	\put(30,33){Absorbing}
	\put(34,28){layers}
	\put(0,40){\line(1,0){50}}
	\put(52,39){$y = \delta$}
	\multiput(0,40)(5,0){10}{\line(1,1){5}}
	\put(25,40){\vector(0,1){10}}
	\put(24,52){$y$}
   \end{picture}
\end{center}

\begin{proposition}[System of acoustic absorbing layers]
\label{pml_1}
\end{proposition}

\nototo
A system of partial differential equations that models absorbing layers of acoustic 
waves in the domain $\Omega^{+}$ introduced in (\ref{def_omega_plus}) can be given as~:

\moneq
\left \lbrace
  \begin{array}{lcl}
        \displaystyle{ \frac{\partial p_x}{\partial t} } & + & 
\displaystyle{ c_0^2 \frac{\partial \xi}{\partial x} = 0 } \\
        \displaystyle{ \frac{\partial p_y}{\partial t} } & + & 
\displaystyle{ \sigma^* (\eta)\, p_y + c_0^2 \frac{\partial \zeta}{\partial \eta} = 0 } \\
        \displaystyle{ \frac{\partial \xi}{\partial t} } & + & 
\displaystyle{ \frac{\partial}{\partial x} (p_x + p_y) = 0 } \\
        \displaystyle{ \frac{\partial \zeta}{\partial t} } & + & 
\displaystyle{ \sigma^* (\eta)\, \zeta + \frac{\partial}{\partial \eta} (p_x + p_y) = 0 }~,
  \end{array}
\right.
\label{sysclabsy}
\end{equation}

\noindent where the absorbing coefficient $\sigma^* (y)$ satisfies~:

\moneq
   [ 0, \delta ] \ni y \longmapsto \sigma^* (y) \in {\textbf{R}}_{+} 
\, \, , \, \,
   \sigma^* (y) > 0 \, \, \, \, \mbox{if} \, \, \, \, y>0 \, \, \, \, 
\mbox{and} \, \, \, \, \sigma^* (0) = 0~.
\label{defsigma}
\end{equation}

The pressure $p$ is defined by $p \equiv p_x + p_y$~.

\begin{proof}
\noindent We follow essentially the idea of Collino \cite{collino}. 
The main idea to establish the acoustic system inside the absorbing layers 
is to introduce the Fourier-Laplace transform and to write the system 
(\ref{sys2_10}) in the complex plan.

\noindent We define the Fourier-Laplace transform by~:

\[
        \displaystyle{ \hat{v}(k_x,y,\omega) = \int \int v(x,y,t) \, 
e^{-i(\omega t + k_x x)} \, {\rm d}x \, {\rm d}t }~.
\]
Then in the domain $\Omega^{-} = \{ (x,y),\mbox{~}y<0 \}$, the system 
(\ref{sys2_10}) takes the form~:

\moneq
\left \lbrace
  \begin{array}{l}
        \displaystyle{ i \omega \hat{p} = -i k_x c_0^2 \hat{\xi} - 
c_0^2 \frac{\partial \hat{\zeta}}{\partial y} } \\
        \displaystyle{ i \omega \hat{\xi} = - i k_x \hat{p} } \\
        \displaystyle{ i \omega \hat{\zeta} = - \frac{\partial \hat{p}}
{\partial y} }~,
  \end{array}
\right.
\label{sys3_2}
\end{equation}

\noindent and the solution of the system (\ref{sys3_2}) for a propagation 
in the growing $y$-direction is obtained after the integration of an ordinary 
differential equation of degree 2~:

\moneq
\left \lbrace
  \begin{array}{l}
        \displaystyle{ \hat{p} = p_0 e^{-i k_y y} } \\
        \displaystyle{ \hat{\xi} = \frac{p_0}{\omega} k_x e^{-i k_y y} } 
\mbox{~~~~~~~~~~~~~~with~~~~}\displaystyle{k_y^2 = \frac{\omega^2}{c_0^2} - k_x^2 }~. \\
        \displaystyle{ \hat{\zeta} = - \frac{p_0}{\omega} k_y e^{-i k_y y} } \\
  \end{array}
\right.
\label{syssol_y_croissant}
\end{equation}

\noindent 
We establish a modified form of the system (\ref{sys3_2}) in $\Omega^{+}$ 
that ensures that waves leaving the domain are not reflected back. Let us extend 
the variable $y$ in the complex plan by adding an imaginary part depending on the 
function $\sigma^{*}$ defined by equation (\ref{defsigma}), and that equals zero 
for $\sigma^{*} \equiv 0$. We write precisely in $\Omega^{+}$ the complex variable 
$y$ parameterized by a real variable $\eta$, $0 \le \eta \le \delta$, as~:
\[
        \displaystyle{ y = \varphi(\eta) = \eta + \frac{1}{i \omega} \int_{0}^{\eta} 
\sigma^* (u) \, {\rm d}u }~.
\]

\noindent We can draw the function $y = \varphi(\eta)$ in the complex plan as~:

\begin{center}
  \begin{picture}(45,29)
	\put(0,20){\line(1,0){40}}
	\put(20,0){\vector(0,1){25}}
	\put(18,27){Re~$y$}
	\put(42,19){$\delta$}
	\put(45,2){\vector(-1,0){37}}
	\put(0,2){Im~$y$}
	\bezier{150}(20,2)(24,16)(32,20)
	\put(16,12){$\eta$}
	\put(19,13){\vector(1,0){5}}
	\put(28,12){$y=\varphi(\eta)$}
  \end{picture}
\end{center}

\noindent For $v$ equal to one of the variables of pressure or momentum, we 
introduce the function $\hat{v} (y)$, with $y = \varphi(\eta)$, as a function 
$\hat{v}^{*}$ of the real variable $\eta$~:

\[
\hat{v}^{*} (\eta) \equiv \hat{v} (\varphi(\eta)) \, \, , \, \, v \in \{ p \, , \xi \, , \zeta \}~.
\]

\noindent An elementary calculation gives us~:

\[
\begin{array}{lcl}
   \hat{v}^{*} (\eta) & = & \displaystyle{ \hat{v} \left( \eta + \frac{1}{i \omega} 
\int_{0}^{\eta} \sigma^* (u) \, {\rm{d}}u \right) } \\
    & = & \displaystyle{ v_0 \, \, \exp \left( -i k_y \left[ \eta + \frac{1}{i 
\omega} \int_{0}^{\eta} \sigma^* (u) \, {\rm{d}}u \right] \right) } \\
    & = & \displaystyle{ \hat{v} \, (\eta) \, \, \exp \left( - \frac{k_y}{\omega} 
\int_{0}^{\eta} \sigma^* (u) \, {\rm{d}}u \right) }~.
\end{array}
\]

\noindent We then deduce the following important property~:

\moneq
	|\hat{v}^* (\eta)| < |\hat{v} \, (\eta)| \, \, \, \, , \, \, \, \, 
0 < \eta \le \delta \, \, \, \, , \, \, \, \, v \in \{ p \, , \xi \, , \zeta \}~.
\label{decayp}
\end{equation}

\noindent The property (\ref{decayp}) is a consequence of an exponential 
decay of all the variables inside the absorbing layers. It is possible to 
derive the system of partial differential equations satisfied by those fields. 
A first algebraic calculation gives us~:

\[ 
\displaystyle{ \frac{\partial \hat{v}}{\partial y} = \frac{\partial \hat{v}^*}
{\partial \eta} \frac{\rm{d} \eta}{{\rm{d}} y} = \frac{i \omega}{i \omega + 
\sigma^*(\eta)} \frac{\partial \hat{v}^*}{\partial \eta} } \, \, \, \, , 
\, \, \, \, v \in \{ p \, , \xi \, , \zeta \}~,
\]

\noindent then we obtain with this new set of unknown functions a system 
in the $(x,\eta)$ domain issued from (\ref{sys3_2})~:

\moneq
\left \lbrace
  \begin{array}{l}
        \displaystyle{ i \omega \hat{p}^* = -i k_x c_0^2 \hat{\xi}^* - c_0^2
\frac{i \omega}{i \omega + \sigma^*(\eta)} \frac{\partial \hat{\zeta}^*}{\partial \eta} } \\
        \displaystyle{ i \omega \hat{\xi}^* = - i k_x \hat{p}^* } \\
        \displaystyle{ i \omega \hat{\zeta}^* = - \frac{i \omega}{i \omega + 
\sigma^*(\eta)} \frac{\partial \hat{p}^*}{\partial \eta} }~.
  \end{array}
\right.
\end{equation}

\noindent The first equation can be rewritten while splitting the pressure field 
into two sub-pressure fields as~:
\[ \hat{p}^* = \hat{p}_x^* + \hat{p}_y^*~, \]
\noindent with $\hat{p}_x^*$ and $\hat{p}_y^*$ solutions of~:
\[
\left \lbrace
   \begin{array}{l}
	\displaystyle{ i \omega \hat{p}_x^* = - i k_x c_0^2 \hat{\xi}^* } \\
	\displaystyle{ i \omega \hat{p}_y^* = \frac{i \omega}{i \omega + 
\sigma^*(\eta)} c_0^2 \frac{\partial \hat{\zeta}^*}{\partial \eta} }~.
   \end{array}
\right.
\]

\noindent Taking the inverse Fourier-Laplace transform of the new 
system, we obtain (\ref{sysclabsy}).
\end{proof}

\noindent
Considering a square domain $\, [ \, 0 \, , \, L \, ]^2$, we define 
the thickness of the absorbing layers by $\delta_x$ in the $x$-direction 
and  $\delta_y$ in the $y$-direction. The interesting studying medium is 
then $\, [ \, \delta_x , L - \delta_x \, ] \, \times \, [ \, \delta_y , L - 
\delta_y \, ]$. We have the following proposition that generalizes the 
proposition \ref{pml_1}~:

\newpage
\begin{proposition}[General acoustic system to solve]
\end{proposition}

\nototo
We consider two smoothing functions $\sigma^{*}_{x}$ and $\sigma^{*}_{y}$ 
defined by~:

\moneq
\begin{array}{l}
    [ 0 , \, L ] \ni x \longmapsto \sigma^{*}_{x} (x) \in {\textbf{R}}_{+}  
	\begin{array}{l}
	     \sigma^{*}_{x} (x) > 0 \,  \mbox{if} \, 
x \in [ \, 0 \, , \delta_x \, [ \, \times \, ] \, L - \delta_x , \, L \, ] \\
	     \sigma^{*}_{x} (x) = 0 \, \,  \mbox{if} \, 
 x \in [ \delta_x , L - \delta_x \, ],
	\end{array}
\end{array}
\end{equation}

\moneq
\begin{array}{l}
    [ 0 , \, L ] \ni \ni y \longmapsto \sigma^{*}_{y} (y) \in {\textbf{R}}_{+}
	\begin{array}{l}
	     \sigma^{*}_{y} (y) > 0 \, \mbox{if} \, 
y \in [ \, 0 \, , \delta_y \, [ \, \times \, ] \, L - \delta_y , \, L \, ] \\
	     \sigma^{*}_{y} (y) = 0 \, \mbox{if} \, 
y \in [ \delta_y , L - \delta_y \, ].
	\end{array}
\end{array}
\label{sigx_sigy}
\end{equation}

\noindent 
The acoustic system in the studying medium and in the absorbing layers can be written as~:

\moneq
  \left \lbrace
    \begin{array}{lclclcl}
        \displaystyle{ \frac{\partial p_x}{\partial t} } & + & 
\displaystyle{ \sigma^{*}_{x} (x) \, p_x } & + & \displaystyle{ c_0^2 
\frac{\partial \xi}{\partial x} } & = & 0 \\
        \displaystyle{ \frac{\partial p_y}{\partial t} } & + & 
\displaystyle{ \sigma^{*}_{y} (y) \, p_y } & + & \displaystyle{ c_0^2  
\frac{\partial \zeta}{\partial y} } & = & 0 \\
        \displaystyle{ \frac{\partial \xi}{\partial t} } & + 
&\displaystyle{ \sigma^{*}_{x} (x) \, \xi } & + & \displaystyle{ 
\frac{\partial}{\partial x} (p_x + p_y) } & =  & 0 \\
        \displaystyle{ \frac{\partial \zeta}{\partial t} } & + 
& \displaystyle{ \sigma^{*}_{y} (y) \, \zeta } & + & \displaystyle{  
\frac{\partial}{\partial y} (p_x + p_y) } & = & 0~.
    \end{array}
  \right.
\label{sys3_5}
\end{equation}

\smallskip \noindent {\bf   Remark 9. } 

\noindent 
This set of equations is the same as obtained by Hu in \cite{hu} using 
velocity fields rather than impulses.

\begin{proof} The proof is similar to the one of proposition \ref{pml_1}
\end{proof}

\bigskip   \noindent {\bf \large 5-2  \quad Plane wave analysis}  

\begin{center}
    \begin{picture}(100,40)
	\put(10,20){\vector(1,0){76}}
	\put(87,19){$x$}
	\put(10,20){\vector(0,1){15}}
	\put(10,36){$y$}
	\put(10,5){\line(2,1){70}}
	\put(40,20){\line(2,-1){40}}
	\put(40,31){2}
	\put(40.8,32){\circle{5}}
	\put(40,7){1}
	\put(40.8,8){\circle{5}}
	\put(10,5){\vector(2,1){10}}
	\put(10,5){\circle*{1}}
	\put(12,2){$(k^{i}_x,k^{i}_y)$}
	\put(70,5){\vector(2,-1){10}}
	\put(70,5){\circle*{1}}
	\put(72,6){$(k^{r}_x,k^{r}_y)$}
	\put(70,35){\vector(2,1){10}}
	\put(70,35){\circle*{1}}
	\put(72,32){$(k^{t}_x,k^{t}_y)$}
	\bezier{50}(47,20)(46.7,22.8)(46,23)
	\put(48,21){$\theta_2$}
	\bezier{70}(49,20)(48.7,16.3)(48,16)
	\put(50,16){$\theta_1'$}
	\bezier{70}(31,20)(31.3,16.3)(32,16)
	\put(28,16.5){$\theta_1$}
    \end{picture}
\end{center}

The acoustic system, inside the absorbing layers in the $y$-direction is, 
after a Fourier transform~:

\[
  \left \lbrace
   \begin{array}{l}
	\displaystyle{ i \omega p_x - i k_x c_0^2 \xi = 0 } \\
	\displaystyle{ i \omega p_y + \sigma^{*} (y) p_y - i k_y c_0^2 \zeta = 0 } \\
	\displaystyle{ i \omega \xi - i k_x p = 0 } \\
	\displaystyle{ i \omega \zeta + \sigma^{*} (y) \zeta - i k_y p = 0 }~.
   \end{array}
  \right.
\]

\noindent The incident wave is a solution of the system~:

\moneq
  \left \lbrace
   \begin{array}{l}
	\displaystyle{ i \omega p - i k^{i}_x c_0^2 \xi - 
\frac{i \omega}{i \omega + \sigma^{*}_1} i k^{i}_y c_0^2 \zeta = 0 } \\
	\displaystyle{ i \omega \xi - i k^{i}_x p = 0 } \\
	\displaystyle{ i \omega \zeta - \frac{i \omega}{i \omega + 
\sigma^{*}_1} i k^{i}_y p = 0 }~.
   \end{array}
  \right.
\label{sys_incident}
\end{equation}

\noindent The reflected wave is a solution of the system~:

\moneq
  \left \lbrace
   \begin{array}{l}
	\displaystyle{ i \omega p_r - i k^{r}_x c_0^2 \xi_r - 
\frac{i \omega}{i \omega + \sigma^{*}_1} i k^{r}_y c_0^2 \zeta_r = 0 } \\
	\displaystyle{ i \omega \xi_r - i k^{r}_x p_r = 0 } \\
	\displaystyle{ i \omega \zeta_r - \frac{i \omega}{i \omega + 
\sigma^{*}_1} i k^{r}_y p_r = 0 }~.
   \end{array}
  \right.
\label{sys_reflechi}
\end{equation}

\noindent The transmitted wave is a solution of the system~:

\moneq
  \left \lbrace
   \begin{array}{l}
	\displaystyle{ i \omega p_t - i k^{t}_x c_0^2 \xi_t - 
\frac{i \omega}{i \omega + \sigma^{*}_2} i k^{t}_y c_0^2 \zeta_t = 0 } \\
	\displaystyle{ i \omega \xi_t - i k^{r}_x p_t = 0 } \\
	\displaystyle{ i \omega \zeta_t - \frac{i \omega}{i \omega +
 \sigma^{*}_2} i k^{r}_y p_t = 0 }~.
   \end{array}
  \right.
\label{sys_transmis}
\end{equation}

\noindent At $y=0$, we write the continuity of the pressure field. We have~:
$ \,\,  	p + p_r = p_t \, \,  $  then  $ \, \, p_r = R \, p \, $ 
and $ \, \, p_t = T p \, $ with   $ \, 1 + R = T .$  
We then have~:

\[
  \left \lbrace
    \begin{array}{l}
	\displaystyle{ \xi_r = R \, \xi \, \, \, \, \mbox{and} \, \, 
\, \, \xi_t = T \xi } \\
	\displaystyle{ \zeta_r = R \, \zeta \, \, \, \, \mbox{and} \, 
\, \, \, \zeta_t = T \zeta }
    \end{array}
  \right.
\]

\noindent We deduce from the system (\ref{sys_incident})~:

\[
	\displaystyle{ \xi = \frac{k^{i}_x}{\omega} p \, \, \, 
\, \, \, \mbox{and} \, \, \, \, \, \, \zeta = \frac{i k^{i}_y}{i 
\omega + \sigma^{*}_1} p }~,
\]
\noindent and we know that $k^{r}_x = k^{i}_x$ and $k^{r}_y = - 
k^{i}_y$. We then have~:

\[
  \left \lbrace
    \begin{array}{l}
	\displaystyle{ \xi + \xi_r = \frac{k^{i}_x}{\omega} p + 
\frac{k^{i}_x}{\omega} p_r = \frac{k^{i}_x}{\omega} (1 + R) p } \\
	\displaystyle{ \zeta + \zeta_r = \frac{i k^{i}_y}{i \omega + 
\sigma^{*}_1} (1 - R) p }
    \end{array}
  \right.
\]

\noindent The system (\ref{sys_transmis}) gives us~:

\[
  \left \lbrace
    \begin{array}{l}
	\displaystyle{ \xi_t = \frac{k^{t}_x}{\omega} (1+R) p = 
\frac{k^{t}_x}{\omega} T p } \\
	\displaystyle{ \xi_t = \frac{i k^{t}_y}{i \omega + \sigma^{*}_2} (1+R) p
 = \frac{i k^{t}_y}{i \omega + \sigma^{*}_2} T p }
    \end{array}
  \right.
\]

\noindent At the interface $y=0$, we write the continuity of $\zeta$. 
We then have $\zeta_t = \zeta + \zeta_r$, that we write as~:

\moneq
	\displaystyle{ \frac{i k^{t}_y}{i \omega + \sigma^{*}_2 (0)} (1 + R) p 
= \frac{i k^{i}_y}{i \omega + \sigma^{*}_1 (0)} (1 - R) p }
\label{tegaliplusr}
\end{equation}

\bigskip  \noindent {\bf \large 5-3  \quad Mathematical property of the absorbing  layers system} 

\noindent 
 The system of acoustic absorbing layers (\ref{sys3_5}) can be written as~:

\moneq
        \displaystyle{ \frac{\partial W}{\partial t} + A 
\frac{\partial W}{\partial x} + B \frac{\partial W}{\partial y} + C \, W = 0 }~,
\label{pml_matrice}
\end{equation}

\noindent where $W = ( p_x,\,p_y,\,\xi,\,\zeta )^t~$, 

\[
A = \left(
          \begin{array}{cccc}
                0 & 0 & c_0^2 & 0 \\
                0 & 0 & 0 & 0 \\
                1 & 1 & 0 & 0 \\
                0 & 0 & 0 & 0 
          \end{array}
      \right),~
B = \left(
          \begin{array}{cccc}
                0 & 0 & 0 & 0 \\
                0 & 0 & 0 & c_0^2 \\
                0 & 0 & 0 & 0 \\
                1 & 1 & 0 & 0 
          \end{array}
       \right),
\]

\noindent  and

\[
C = \left(
          \begin{array}{cccc}
                \sigma_x^*(x) & 0 & 0 & 0 \\
                0 & \sigma_y^*(y) & 0 & 0 \\
                0 & 0 & \sigma_x^*(x) & 0 \\
                0 & 0 & 0 & \sigma_y^*(y)
          \end{array}
      \right).
\]

\noindent The principal symbol, $M = ik_x A + i k_y B$, of the system 
(\ref{pml_matrice}) is given by~:

\moneq
	M = 
   \left(
      \begin{array}{cccc}
        0 & 0 & i k_x c_0^2 & 0 \\
        0 & 0 & 0 & i k_y c_0^2 \\
        i k_x & i k_x & 0 & 0 \\
        i k_y & i k_y & 0 & 0
      \end{array}
   \right)~.
\label{eigenvalues}
\end{equation}

\noindent The eigenvalues and eigenvectors of the principal symbol are
  re-evaluated       in table 1. 

\begin{tabular*}{\textwidth}{@{\extracolsep{\fill}}cccc}  
\\ \hline \\  
 [-5mm]   & Eigenvalue & Eigenvector & 
\\[1mm]\hline \\[-1mm] 
 & 0 (double) & $(-1,1,0,0)^t$ & \\  [2mm] 
 & $ \displaystyle c_0 \sqrt{-k_x^2 - k_y^2}$ & $ \displaystyle ( -\frac{i k_x c_0 \sqrt{-k_x^2 - 
k_y^2}}{k_x^2 + k_y^2},-\frac{i k_y^2 c_0 \sqrt{-k_x^2 - k_y^2}}{ kx 
( k_x^2 + k_y^2 ) },1,\frac{k_y}{k_x} )^t$ & \\ [3mm] 
 & $ \displaystyle - c_0 \sqrt{-k_x^2 - k_y^2}$ & $  \displaystyle ( \frac{i k_x c_0 \sqrt{-k_x^2 - 
k_y^2}}{k_x^2 + k_y^2},\frac{i k_y^2 c_0 \sqrt{-k_x^2 - k_y^2}}{k_x 
( k_x^2 + k_y^2 ) },1,\frac{k_y}{k_x} )^t$ & 
\\[2mm]  \hline 
\end{tabular*}
~  \bigskip  \centerline { {\bf Table 1}. \quad  
Eigenvalues and eigenvectors of matrix (\ref{eigenvalues})}
 \bigskip 

\noindent 
We notice that 0 is an eigenvalue of multiplicity order equal to 2 
associated with a one-dimensional eigensubspace. The system (\ref{pml_matrice}) 
is not hyperbolic and classical results relative to well-posedness of such 
systems (see \cite{rauch}) can not be applied to the absorbing layers. 
There exists an arbitrarily small perturbation of the Cauchy problem for 
the system (\ref{sys3_5}) with $\sigma^*_{x} (x) = \sigma^*_{y} (y) = 0$ 
that makes the system (\ref{pml_matrice}) ill-posed for $L^2$ or Sobolev 
norms of order 1. Nevertheless, our choice of the system (\ref{pml_matrice}) 
does not produce unstable numerical results.

\begin{proposition}
\end{proposition}

\nototo
If we look for a solution of the form $W = \varphi (t) \, e^{- i k_x x}
 e^{- i k_y y} \, V_{M^2}$ of the system (\ref{pml_matrice}), where $V_{M^2} = 
(0, 0, k_y, -k_x)^{t}$, the scalar function $\varphi (t)$ is an exponential decay in time.

\begin{proof}
To establish this result, we determine the characteristic subspace 
for the eigenvalue $\lambda = 0$, {\it{i.e.}} we calculate $ker(M^2)$. We have~:
	
\moneq
- M^2 = 
   \left(
      \begin{array}{cccc}
        {k_x}^2 & {k_x}^2 & 0 & 0 \\
        {k_y}^2 & {k_y}^2 & 0 & 0 \\
        0 & 0 & {k_x}^2 & k_x k_y \\
        0 & 0 & k_x k_y & {k_y}^2
      \end{array}
   \right) \,  \mbox{and} \quad  {\rm ker} (M^2) = 
	\left[
		\left(
			\begin{array}{c}
				1 \\ -1 \\ 0 \\ 0
			\end{array}
		\right) ~,~
		\left(
			\begin{array}{c}
				0 \\ 0 \\ k_y \\ - k_x
			\end{array}
		\right)
	\right]~.
\end{equation}

\noindent We analyse the stability of the system (\ref{pml_matrice}) 
under a perturbation following the direction of the eigenvector of $M^2$ 
that is not eigenvector of $M$. We note $V_{M^2} = (0, 0, k_y, -k_x)^{t}$ 
this vector that is a simple impulse pertubation in the direction orthogonal 
to the wave vector. We choose a state vector $W$ of the form~:
\[ W = \varphi (t) \, e^{- i k_x x} e^{- i k_y y} \, V_{M^2}~. \]
\noindent In this case, the system (\ref{sys3_5}) is written as~:

\moneq 
	\frac{\partial \varphi}{\partial t} \, V_{M^2} + \varphi M \, V_{M^2} 
+ \varphi C \, V_{M^2} = 0~,
\label{pertmat}
\end{equation}

\noindent with~:

\[ 
   MV_{M^2} = 
	\left(
		\begin{array}{c}
			i k_x k_y \\ - i k_x k_y \\ 0 \\ 0
		\end{array}
	\right) ~, ~
   CV_{M^2} = 
	\left(
		\begin{array}{c}
			0 \\ 0 \\ \sigma^{*}_x (x) k_y \\ - \sigma^{*}_y (y) k_x
		\end{array}
	\right)~.
\]

\noindent We then deduce that the first two equations of (\ref{pertmat}) 
impose that $k_x k_y = 0$. Therefore, if $k_x = 0$ and $k_y \neq 0$, 
the third equation of (\ref{pertmat}) gives us~:

\moneq
	\frac{\partial \varphi}{\partial t} k_y + \sigma^{*}_x (x) k_y 
\varphi = 0 \, \, \, \, \mbox{{\it{i.e.}}} \, \, \, \, \frac{\partial 
\varphi}{\partial t} + \sigma^{*}_x (x) \varphi = 0~.
\label{eq_dis_pml}
\end{equation}

\noindent We assume that $\sigma^{*}_x (x) > 0$ inside the absorbing 
layers, the solution of the equation (\ref{eq_dis_pml}) is an exponential 
decay in time of the function $\varphi$.
\end{proof}

\smallskip \noindent 
We obtain the same result considering $k_y = 0$ and $k_x \neq 0$. 
Even if the principal symbol of the system (\ref{pml_matrice}) is 
associated to a ``ill-posed mathematical problem'', the form of the 
zero order terms shows that even exciting the absorbing layers system in 
the direction of the characteristic vector, the perturbation is dissipated.

\smallskip \noindent 
We would like to predict the behavior of our absorbing layers model. 
Thus, we simplify our set of equations to the simplest model and study it. 
We establish the following proposition~:

\begin{proposition}
\label{prop_1D}
\end{proposition}

 \nototo
We consider the simplest 1-D non-hyperbolic model inside the absorbing layers, 
excitated with a source term $\psi (t)$ centered at $(x_a, y_a)$, in the 
direction of the eigenvector $V_{M^2}$. We note $W = (u,v)^t$ the state vector, 
$\delta_{x_a,y_a}$ the Dirac mass at the position $(x_a, y_a)$ and we assume 
that the coefficients $\sigma_1$ and $\sigma_2$ are strictly positive. The problem~:

\moneq 
  \left \lbrace
    \begin{array}{l}
        \displaystyle{ \frac{\partial W}{\partial t} + 
\left( \begin{array}{cc}
        0 & 1 \\ 0 & 0
\end{array} \right) \frac{\partial W}{\partial x} + 
\left( \begin{array}{cc}
        \sigma_1 & 0 \\ 0 & \sigma_2
\end{array} \right) W = 
\left( \begin{array}{c}
        0 \\ \psi(t) \, \delta_{x_a,y_a}
\end{array} \right)
 } \\
        W (0) = 0~.
    \end{array}
  \right.
\label{sol_pml}
\end{equation}
is stable as long as $\psi (t)$ is bounded.

\begin{proof}
The system (\ref{sol_pml}) is written as~:

\moneq
  \left \lbrace
    \begin{array}{l}
        \displaystyle{ \frac{\partial u}{\partial t} + 
\frac{\partial v}{\partial x} + \sigma_1 u = 0 } \\
        \displaystyle{ \frac{\partial v}{\partial t} + 
\sigma_2 v = \psi(t) \, \delta_{x_a,y_a} } \\
        u (0) = v (0) = 0~.
    \end{array}
  \right.
\label{proof_pml}
\end{equation}

\noindent The solution in $v$ is $\displaystyle{ v (t) = \left( \int_0^t \psi 
(\theta) e^{-\sigma_2 (t - \theta)} {\rm d} \theta \right) \delta_{x_a,y_a} }$, 
we then deduce that $u$ is solution of~:

\[
        \displaystyle{ \frac{\partial u}{\partial t} + \sigma_1 u = 
\left( \int_0^t \psi (\theta) e^{-\sigma_2 (t - \theta)} {\rm d} 
\theta \right) \delta_{x_a,y_a}' }~.
\]

\noindent Then, the solution in $u$ is $u(t) = \mu (t) \, \delta_{x_a,y_a}'$, 
where $\lim_{t \to \infty} \mu(t) = 0$. Then, for a function $\psi (t)$ 
bounded, $u \to 0$ and $v$ is bounded. The system (\ref{proof_pml}) is
 stable as long as $\sigma_1 > 0$ and $\sigma_2 > 0$.
\end{proof}

\smallskip \noindent 
This simple model makes us think that our absorbing layers model is 
stable even if we are in the worst situation, {\it{i.e.}} there is a 
source in the direction of the eigenvector $V_{M^2}$. Numerical simulations 
proposed in Section 7 
confirm this result.

\smallskip \noindent 
One challenge in the future is to theorically understand the behavior of 
the solution of the true problem in the absorbing layers defined by the following system~:

\moneq
  \left \lbrace
    \begin{array}{l}
        \displaystyle{ \frac{\partial W}{\partial t} + A \frac{\partial W}
{\partial x} + B \frac{\partial W}{\partial y} + C \, W = \psi(t) \, 
\delta_{x_a,y_a} \, V_{M^2} } \\
        W (0) = 0~.
    \end{array}
  \right.
\label{truePb}
\end{equation}

\noindent Nevertheless, the qualitative behavior proposed for the system 
(\ref{sol_pml}) gives a good idea of the behavior of the system (\ref{truePb}), 
see Section 7.  

\bigskip  \noindent {\bf \large 5-4  \quad Acoustic absorbing layers with external subsonic flow}

\noindent 
\begin{proposition}[General absorbing layers in two dimensions]
\end{proposition}

\nototo
We assume that the velocity for the external subsonic flow is ${\bf{u}} = (u_0, v_0)$~. 
A system of partial differential equations that models absorbing layers of 
acoustic waves is given by~:

\moneq
    \left \lbrace
        \begin{array}{l}
          \displaystyle{ \frac{\partial p_x}{\partial t} + c_0 \sqrt{1-M_0^2} \, 
\sigma^* (x) \, p_x + \frac{c_0 u_0}{\sqrt{1-M_0^2}} \sigma^* (x) \, \xi + c_0^2 
\frac{\partial \xi}{\partial x} = 0 } \\
          \displaystyle{ \frac{\partial p_y}{\partial t} + c_0 \sqrt{1-M_0^2} \, 
\sigma^* (y) \, p_y + \frac{c_0 v_0}{\sqrt{1-M_0^2}} \sigma^* (y) \, \zeta + c_0^2
 \frac{\partial \zeta}{\partial y} = 0 } \\
          \begin{array}{lcl}
                \displaystyle{ \frac{\partial \xi}{\partial t} } & + & \displaystyle{ 
\frac{\partial}{\partial x} ( 2 u_0 \xi + v_0 \zeta) + (1-M_0^2) \frac{\partial p}
{\partial x} + \frac{\partial}{\partial y} (u_0 \zeta) + \frac{c_0 (1+\frac{u_0^2 -
 v_0^2}{c_0^2})}{\sqrt{1-M_0^2}} \sigma^* (x) \xi } \\
                 & + & \displaystyle {\frac{u_0 \sqrt{1-M_0^2}}{c_0} \left( 
\sigma^* (x) p_x + \sigma^* (y) p_y \right) + \frac{u_0 v_0 ( \sigma^* (x) + 
\sigma^* (y) )}{c_0 \sqrt{1-M_0^2}} \zeta = 0 }
          \end{array} \\
          \begin{array}{lcl}
                \displaystyle{ \frac{\partial \zeta}{\partial t} } & + & \displaystyle{ 
\frac{\partial}{\partial x} (v_0 \xi) + \frac{\partial}{\partial y} (u_0 \xi + 2 v_0
\zeta) 
+ (1-M_0^2) \frac{\partial p}{\partial y} + \frac{c_0 (1+\frac{v_0^2-u_0^2}{c_0^2})}
{\sqrt{1-M_0^2}} \sigma^* (y) \zeta } \\
         & + & \displaystyle{ \frac{v_0 \sqrt{1-M_0^2}}{c_0} \left( \sigma^* (x) \, p_x + 
 \sigma^* (y) \, p_y \right) + \frac{u_0 v_0 ( \sigma^* (x) + \sigma^* (y) )}{c_0 
\sqrt{1-M_0^2}} \xi = 0 }~,
          \end{array}
        \end{array}
    \right.
\label{sys_pml_uv}
\end{equation}

\noindent where $\sigma^* (x)$ et $\sigma^* (y)$ are smoothing functions defined 
by (\ref{sigx_sigy}).

\begin{proof} Here are the main ideas of the proof; the details of all the 
calculus can be found in \cite{fred}. If we consider a two-dimensional flow, we 
have shown that the acoustic system in the new space-time defined by (\ref{sys2_11}) 
is given by (\ref{sys2_12}) after the change of unknown functions (\ref{sys2_14}). 
Using the method described in the section, we easily show that a general system of
 dimensionless partial differential equations for the acoustic absorbing layers can 
be written in $(x',\,y',\,t')$ using the functions $\sigma^* (x)$ and $\sigma^* (y)$ 
introduced by equation (\ref{sigx_sigy}) as~:

\moneq
\left \lbrace
  \begin{array}{lclcl}
        \displaystyle{ \frac{\partial \widetilde{p_x}}{\partial t'} } & + & 
\displaystyle{ c_0 \sqrt{1-M_0^2} \, \sigma^* (x) \widetilde{p_x} } & + & 
\displaystyle{ c_0^2 \frac{\partial}{\partial x'} (\widetilde{\xi} - \alpha
 \widetilde{\zeta}) = 0 } \\
        \displaystyle{ \frac{\partial \widetilde{p_y}}{\partial t'} } & + & 
\displaystyle{ c_0 \sqrt{1-M_0^2} \, \sigma^* (y) \widetilde{p_y} } & + &  
\displaystyle{ c_0^2 \frac{\partial}{\partial y'} (\widetilde{\zeta} - \alpha 
\widetilde{\xi}) = 0 } \\
        \displaystyle{ \frac{\partial \widetilde{\xi}}{\partial t'} } & + & 
\displaystyle{ c_0 \sqrt{1-M_0^2} \, \sigma^* (x) \widetilde{\xi} } & + & 
\displaystyle{ \frac{\partial}{\partial x'} (\widetilde{p_x} + \widetilde{p_y}) = 0 } \\
        \displaystyle{ \frac{\partial \widetilde{\zeta}}{\partial t'} } & + & 
\displaystyle{ c_0 \sqrt{1-M_0^2} \, \sigma^* (y) \widetilde{\zeta} } & + & 
 \displaystyle{ \frac{\partial}{\partial y'} (\widetilde{p_x} + \widetilde{p_y}) = 0 }~,
  \end{array}
\right.
\label{sys4_19_1}
\end{equation}

\noindent where $\alpha$ is a coupling coefficient given by~:

\[ \alpha = \frac{u_0 v_0}{c_0^2 \sqrt{1-\frac{u_0^2}{c_0^2}} \sqrt{1-\frac{v_0^2}{c_0^2}}}~. \]

\noindent We now wish to write the system (\ref{sys4_19_1}) in the initial 
space-time $(x,\,y,\,t)$, using the initial unknown functions $p,\, \xi, \, \zeta$. We
have~:

\moneq
        \displaystyle{ \widetilde{\xi} - \alpha \widetilde{\zeta} = \frac{1}
{\sqrt{1-\frac{u_0^2}{c_0^2}}} \, \xi' } \, \, \, \, , \, \, \, \,
        \displaystyle{ \widetilde{\zeta} - \alpha \widetilde{\xi} = 
\frac{1}{\sqrt{1-\frac{v_0^2}{c_0^2}}} \, \zeta' }~,
\label{sys4_19_2}
\end{equation}

\noindent and we easily calculate~:

\moneq
\left \lbrace
  \begin{array}{l}
        \displaystyle{ \frac{\partial}{\partial t'} = \frac{\partial}{\partial t} }  \\
        \displaystyle{ \frac{\partial}{\partial x'} = \sqrt{1-\frac{u_0^2}{c_0^2}} 
\frac{\partial}{\partial x} - \frac{u_0 \sqrt{1-\frac{u_0^2}{c_0^2}}}{c_0^2 (1-M_0^2)} 
\frac{\partial}{\partial t} }  \\
        \displaystyle{ \frac{\partial}{\partial y'} = \sqrt{1-\frac{v_0^2}{c_0^2}}
 \frac{\partial}{\partial y} - \frac{v_0 \sqrt{1-\frac{v_0^2}{c_0^2}}}{c_0^2 (1-M_0^2)} 
\frac{\partial}{\partial t} }~.
  \end{array}
\right.
\label{sys4_19_3}
\end{equation}

\noindent We substitute (\ref{sys4_19_2}) and (\ref{sys4_19_3}) in the system 
(\ref{sys4_19_1})~:

\[
  \left  \lbrace
    \begin{array}{l}
        \displaystyle{ \frac{\partial \widetilde{p_x}}{\partial t} + c_0 
\sqrt{1-M_0^2} \sigma^* (x) \widetilde{p_x} + c_0^2 \sqrt{1-\frac{u_0^2}{c_0^2}} 
\frac{\partial}{\partial x} (\widetilde{\xi} - \alpha \widetilde{\zeta}) - 
\frac{u_0 \sqrt{1-\frac{u_0^2}{c_0^2}}}{1-M_0^2} \frac{\partial}{\partial t} 
(\widetilde{\xi} - \alpha \widetilde{\zeta}) = 0 } \\
        \displaystyle{ \frac{\partial \widetilde{p_y}}{\partial t} + c_0 
\sqrt{1-M_0^2} \sigma^* (y) \widetilde{p_y} + c_0^2 \sqrt{1-\frac{v_0^2}{c_0^2}} 
\frac{\partial}{\partial y} (\widetilde{\zeta} - \alpha \widetilde{\xi}) - 
\frac{v_0 \sqrt{1-\frac{v_0^2}{c_0^2}}}{1-M_0^2} \frac{\partial}{\partial t}
 (\widetilde{\zeta} - \alpha \widetilde{\xi}) = 0 } \\
        \displaystyle{ \frac{\partial \widetilde{\xi}}{\partial t} + c_0 
\sqrt{1-M_0^2} \sigma^* (x) \widetilde{\xi} + (1-M_0^2) \sqrt{1-\frac{u_0^2}{c_0^2}} 
\frac{\partial \widetilde{p}}{\partial x} - \frac{u_0 \sqrt{1-\frac{u_0^2}{c_0^2}}}{c_0^2} 
\frac{\partial \widetilde{p}}{\partial t} = 0 } \\
        \displaystyle{ \frac{\partial \widetilde{\zeta}}{\partial t} + c_0 
\sqrt{1-M_0^2} \sigma^* (y) \widetilde{\zeta} + (1-M_0^2) \sqrt{1-\frac{v_0^2}{c_0^2}} 
\frac{\partial \widetilde{p}}{\partial y} - \frac{v_0 \sqrt{1-\frac{v_0^2}{c_0^2}}}{c_0^2} 
\frac{\partial \widetilde{p}}{\partial t} = 0 }~.
    \end{array}
  \right.
\]

\noindent Using the change of variables (\ref{sys2_14}), we deduce~:

\moneq
\left \lbrace
  \begin{array}{l}
        \displaystyle{ \widetilde{p} = p' + \frac{1}{1-M_0^2} \left( u_0 \xi' 
+ v_0 \zeta' \right) } \\ 
        \displaystyle{ \widetilde{\xi} = \sqrt{ 1-\frac{u_0^2}{c_0^2} }
 \left( (1 - \frac{v_0^2}{c_0^2}) \frac{\xi'}{1-M_0^2} + \frac{u_0 v_0}{c_0^2} 
\frac{\zeta'}{1-M_0^2} \right) } \\
        \displaystyle{ \widetilde{\zeta} = \sqrt{ 1-\frac{v_0^2}{c_0^2} }
 \left( \frac{u_0 v_0}{c_0^2} \frac{\xi'}{1-M_0^2} + (1 - \frac{u_0^2}{c_0^2}) 
\frac{\zeta'}{1-M_0^2} \right) }~.
  \end{array}
\right.
\label{change_of_variables}
\end{equation}

\noindent We substitute the change of variables (\ref{change_of_variables}) 
into the previous system, we then obtain~:

\moneq
  \left  \lbrace
    \begin{array}{l}
        \displaystyle{ \frac{\partial \widetilde{p_x}}{\partial t} - 
\frac{u_0}{1-M_0^2} \frac{\partial \xi}{\partial t} + c_0 \sqrt{1-M_0^2} 
\sigma^* (x) \widetilde{p_x} + c_0^2 \frac{\partial \xi}{\partial x} = 0 } \\
        \displaystyle{ \frac{\partial \widetilde{p_y}}{\partial t} - 
\frac{v_0}{1-M_0^2} \frac{\partial \zeta}{\partial t} + c_0 \sqrt{1-M_0^2} 
\sigma^* (y) \widetilde{p_y} + c_0^2 \frac{\partial \zeta}{\partial y} = 0 } \\
        \begin{array}{lcl}
	\displaystyle{ \frac{\partial \xi}{\partial t} - \frac{u_0}{c_0^2} 
\frac{\partial p}{\partial t} } & + & \displaystyle{ c_0 \sqrt{ 1-M_0^2 } 
\sigma^* (x) \left( (1-\frac{v_0^2}{c_0^2}) \frac{\xi}{1-M_0^2} +
\frac{u_0 v_0}{c_0^2} \frac{\zeta}{1-M_0^2} \right) } \\
               & + & \displaystyle{ (1-M_0^2) \frac{\partial p}{\partial x} + 
\frac{\partial}{\partial x} (u_0 \xi + v_0 \zeta) = 0 } \\
        \end{array} \\
        \begin{array}{lcl}
             \displaystyle{ \frac{\partial \zeta}{\partial t} - 
\frac{v_0}{c_0^2} \frac{\partial p}{\partial t} } & + & \displaystyle{ c_0 
\sqrt{ 1-M_0^2 } \sigma^* (y) \left( \frac{u_0 v_0}{c_0^2} \frac{\xi}{1-M_0^2} 
+ (1-\frac{u_0^2}{c_0^2}) \frac{\zeta}{1-M_0^2} \right) } \\
               & + & \displaystyle{ \frac{\partial}{\partial y} (u_0 \xi + v_0 
\zeta) + (1-M_0^2) \frac{\partial p}{\partial y } = 0 }~.
        \end{array}
    \end{array}
  \right.
\label{sys4_19_5}
\end{equation}

\noindent We use as new unknowns~:

\[
        \displaystyle{ \widetilde{p_x} = p_x + \frac{u_0}{1-M_0^2} \xi } 
\, \, \, \, , \, \, \, \,
        \displaystyle{ \widetilde{p_y} = p_y + \frac{v_0}{1-M_0^2} \zeta }~,
\]

\noindent we then have~:
\[
        \displaystyle{ \widetilde{p} = \widetilde{p_x} + \widetilde{p_y} } \
, \, \, \, \, \, \, \, \,
        \displaystyle{ p = p_x + p_y }~,
\]

\noindent and we finally obtain~:

\moneq
  \left  \lbrace
    \begin{array}{l}
        \displaystyle{ \frac{\partial p_x}{\partial t} + c_0 \sqrt{1-M_0^2} 
\sigma^* (x) p_x + c_0^2 \frac{\partial \xi}{\partial x} = 0 } \\
        \displaystyle{ \frac{\partial p_y}{\partial t} + c_0 \sqrt{1-M_0^2} 
\sigma^* (y) p_y + \frac{c_0 v_0}{\sqrt{1-M_0^2}} \sigma^* (y) \zeta + c_0^2 
\frac{\partial \zeta}{\partial y} = 0 } \\
        \displaystyle{ \frac{\partial \xi}{\partial t} - \frac{u_0}{c_0^2} 
\frac{\partial p}{\partial t} + (1-M_0^2) \frac{\partial p}{\partial x} + 
\frac{\partial}{\partial x} (u_0 \xi + v_0 \zeta) = 0 } \\
        \begin{array}{lcl}
             \displaystyle{ \frac{\partial \zeta}{\partial t} - \frac{v_0}{c_0^2} 
\frac{\partial p}{\partial t} } & + & \displaystyle{ c_0 \sqrt{ 1-M_0^2 } 
\sigma^* (y) \left( \frac{u_0 v_0}{c_0^2} \frac{\xi}{1-M_0^2} + 
(1-\frac{u_0^2}{c_0^2}) \frac{\zeta}{1-M_0^2} \right) } \\
               & + & \displaystyle{ \frac{\partial}{\partial y} (u_0 \xi + 
v_0 \zeta) + (1-M_0^2) \frac{\partial p}{\partial y } = 0 }~.
        \end{array}
    \end{array}
  \right.
\label{sys4_19_6}
\end{equation}

\noindent We wish to find a dynamic system, then we eliminate the $\frac{\partial 
p}{\partial t}$ terms in the last two equations of the system (\ref{sys4_19_6}). 
Using the equality $p = p_x + p_y$, and adding the first two equations of the 
system (\ref{sys4_19_6}), we deduce~:

\moneq
        \displaystyle{ \frac{\partial p}{\partial t} + c_0^2 \frac{\partial \xi}
{\partial x} + c_0^2 \frac{\partial \zeta}{\partial y} + c_0 \sqrt{1-M_0^2} 
\sigma^* (y) p_y + \frac{c_0 v_0}{\sqrt{1-M_0^2}} \sigma^* (y) \zeta = 0 }~,
\label{sys4_19_7}
\end{equation}

\noindent and we substitute the equation (\ref{sys4_19_7}) in the last 
two equations of the system (\ref{sys4_19_6}). Hence, we find the result
 (\ref{sys_pml_uv}) that ends the proof.
\end{proof}
 
\smallskip \noindent {\bf   Remark 13. } 

\noindent 
Various authors propose a system of partial differential equation for 
absorbing layers for advective acoustic (see \cite{gottlieb}, \cite{hu}, 
\cite{adib1}, \cite{adib2} for example). Each of them have to solve 6 equations, 
whereas we propose a system composed by only 4 equations. We see this 
property as a consequence of our precise physical analysis based on the 
Lorentz transform and our change of unknown functions.

\bigskip \bigskip  \noindent {\bf \large 6) \quad Discretization with the ``HaWAY'' method}

\noindent 
This section deals with the numerical resolution of the equation of 
acoustic {\bf{without}} an external flow. We use finite differences with 
staggered grids as introduced by Harlow-Welsch (MAC method) \cite{harlow}, 
Arakawa (C grids) \cite{arakawa} and Yee \cite{yee} for electromagnetism. 
HaWAY comes from {\bf{Ha}}rlow-{\bf{W}}elsch, {\bf{A}}rakawa, {\bf{Y}}ee. 
The acoustic system can be written as~:

\moneq
\left \lbrace
  \begin{array}{l}
    \displaystyle{ \frac{\partial p'}{\partial t'}} + c_0^2 
\frac{\partial \xi'}{\partial x'} + c_0^2 \frac{\partial \zeta'}{\partial y'} = 0 \\
    \displaystyle{ \frac{\partial \xi'}{\partial t'} + \frac{\partial p'}
{\partial x'} = 0 } \\
    \displaystyle{ \frac{\partial \zeta'}{\partial t'} + \frac{\partial p'}
{\partial y'} = 0 }~.
  \end{array}
\right.
\label{sys4_1}
\end{equation}

\noindent 
We first propose to nondimensionalize the previous system. We then explain 
the numerical scheme chosen in the free space and in the acoustic absorbing layers.

\bigskip   \noindent {\bf \large 6-1 \quad Dimensionlessness of the acoustic system}

\noindent 
This section is introduced for the completeness of our meaning. We refer 
to \cite{sedov} for this kind of purpose. We nondimensionalize the set of 
equations (\ref{sys4_1}) by writing each variable as~:
$ \,\,   X' = X^* X, \, $   for $X'$ pressure, impulse, time and space variables 
and $X^*$ a reference 
dimension for each variable: a reference pressure $p^*$, reference impulses $\xi^*$ 
and $\zeta^*$, a time reference $t^*$ and reference lengths $x^*$ and $y^*$. 
The new form of the system (\ref{sys4_1}) is~:

\[
\left \lbrace
  \begin{array}{l}
        \displaystyle{ \frac{\partial p}{\partial t} + c_0^2 \frac{t^*}{p^*} 
\frac{\xi^*}{x^*} \frac{\partial \xi}{\partial x} + c_0^2 \frac{t^*}{p^*} 
\frac{\zeta^*}{y^*} \frac{\partial \zeta}{\partial y} = 0 } \\
        \displaystyle{ \frac{\partial \xi}{\partial t} + \frac{t^*}{\xi^*} 
\frac{p^*}{x^*} \frac{\partial p}{\partial x} = 0 } \\
        \displaystyle{ \frac{\partial \zeta}{\partial t} + \frac{t^*}{\zeta^*} 
\frac{p^*}{y^*} \frac{\partial p}{\partial y} = 0 }~.
  \end{array}
\right.
\]

\noindent 
We decide here to {\bf{choose}} the following coefficients 
$\displaystyle{ c_0^2 \frac{t^*}{p^*} \frac{\xi^*}{x^*} }$, $\displaystyle{ 
c_0^2 \frac{t^*}{p^*} \frac{\zeta^*}{y^*} }$, $\displaystyle{ \frac{t^*}{\xi^*} 
\frac{p^*}{x^*} }$ and $\displaystyle{ \frac{t^*}{\zeta^*} \frac{p^*}{y^*} }$ 
equal to 1. We then deduce that we have~:

\[
\displaystyle{ \xi^* = \frac{1}{c_0} p^* } \; \; \mbox{,} \; \; 
\displaystyle{ \zeta^* = \frac{1}{c_0} p^* } \; \; \mbox{,} \; \; 
\displaystyle{ \frac{x^*}{t^*} = c_0 } \; \; \mbox{,} \; \; 
\displaystyle{ \frac{y^*}{t^*} = c_0 }~,
\]

\noindent and the resulting set of dimensionless equations is~:

\moneq
\left \lbrace
  \begin{array}{l}
        \displaystyle{ \frac{\partial p}{\partial t} + \frac{\partial \xi}
{\partial x} + \frac{\partial \zeta}{\partial y} = 0 } \\
        \displaystyle{ \frac{\partial \xi}{\partial t} + \frac{\partial p}
{\partial x} = 0 } \\
        \displaystyle{ \frac{\partial \zeta}{\partial t} + \frac{\partial p}
{\partial y} = 0 }~.
  \end{array}
\right.
\label{sys4_2}
\end{equation}

\bigskip   \noindent {\bf \large 6-2 \quad Staggered grids for acoustics }

\noindent 
By analogy with electromagnetism (see \cite{cray}), we decide to use the cartesian 
staggered finite differences method to solve the system (\ref{sys4_2}). We decompose 
a model domain $\Omega = \, [ \, 0 , L \, ]^2$ into finite elements with an isotropic 
meshing of space step ${\Delta x} = {\Delta y} = \frac{L}{J}$ ($J \in {\bf{N^*}}$ is 
the number of cells in each direction). The cell $K_{i+\frac{1}{2},j+\frac{1}{2}}$ is 
defined as~:

\[
	K_{i+\frac{1}{2},j+\frac{1}{2}} = ] \, i {\Delta x} \, , \, (i+1) {\Delta x} 
\, [ \, \times \, ] \, j {\Delta y} \, , \, (j+1) {\Delta y} \, [
\]

\noindent 
We share the time with the help of a time step ${\Delta t}$ and introduce the 
$n^{th}$ ``entire time'' $ \, t^n = n {\Delta t} \, $. By convention, we know that 
the pressure is defined at entire times $t^n$ in the center of the mesh 
$ \, K_{i+\frac{1}{2},j+\frac{1}{2}} \, $ and that the impulses are defined at 
semi-entire times $t^{n+\frac{1}{2}}$ on the edge of the mesh.
 The variables in a mesh are defined as below~:

\centerline { 
  \begin{picture}(60,45)
        \put(8,8){\line(0,1){30}}
        \put(8,8){\line(1,0){30}}
        \put(38,38){\line(0,-1){30}}
        \put(38,38){\line(-1,0){30}}
        \put(23,23){\circle*{2}}
        \put(8,23){\circle*{2}}
        \put(23,8){\circle*{2}}
        \put(22,26){$p^{n}_{i+\frac{1}{2},j+\frac{1}{2}}$}
        \put(10,23){$\xi^{n+\frac{1}{2}}_{i,j+\frac{1}{2}}$}
        \put(22,12){$\zeta^{n+\frac{1}{2}}_{i+\frac{1}{2},j}$}
  \end{picture} }  

\noindent The numerical scheme used in the free space is in two-dimensional space~:

\noindent $\bullet$  \quad   Discretization of the pressure equation~:
   \[ 
        p_{i+ \frac{1}{2},j+ \frac{1}{2}}^{n+1} = p_{i+ \frac{1}{2},j+ 
\frac{1}{2}}^{n} - \sigma \left[ ( \xi_{i+1,j+ \frac{1}{2}}^{n+ \frac{1}{2}} 
- \xi_{i,j+ \frac{1}{2}}^{n+ \frac{1}{2}} ) + ( \zeta_{i+ \frac{1}{2},j+1}^{n+ 
\frac{1}{2}} - \zeta_{i+ \frac{1}{2},j}^{n+ \frac{1}{2}} ) \right]~,
   \]

\noindent $\bullet$  \quad    Discretization of the impulse equations~:
   \[
   \left \lbrace
      \begin{array}{l}
        \xi_{i,j+\frac{1}{2}}^{n+ \frac{3}{2}} = \xi_{i,j+\frac{1}{2}}^
{n+ \frac{1}{2}} - \sigma ( p_{i+\frac{1}{2},j+\frac{1}{2}}^{n+1} - 
p_{i-\frac{1}{2},j+\frac{1}{2}}^{n+1} ) \\ 
        \zeta_{i+\frac{1}{2},j}^{n+ \frac{3}{2}} = \zeta_{i+\frac{1}{2},j}^{n+ 
\frac{1}{2}} - \sigma ( p_{i+\frac{1}{2},j+\frac{1}{2}}^{n+1} - p_{i+\frac{1}{2},
j-\frac{1}{2}}^{n+1} )~,
      \end{array}
   \right.
   \] 

\noindent 
where we have $\displaystyle{\sigma = \frac{{\Delta t}}{{\Delta x}} = 
\frac{{\Delta t}}{{\Delta y}}}$. \\

\noindent 
This numerical scheme is an explicit second order in time and space scheme, stable under the 
Courant-Friedrichs-Lewy condition~:

\moneq
   \begin{array}{ll}
      \mbox{CFL} & \left \lbrace \begin{array}{ll}
        \displaystyle{ {\Delta t} \leq \frac{1}{\sqrt{ \frac{1}{{\Delta x}^2} 
+ \frac{1}{{\Delta y}^2}} } } & \mbox{~~~in two dimension space,} \\
        \displaystyle{ {\Delta t} \leq \frac{1}{\sqrt{\frac{1}{{\Delta x}^2} 
+ \frac{1}{{\Delta y}^2} + \frac{1}{{\Delta z}^2} } } } & \mbox{~~~in three dimension space.}
                          \end{array}
            \right.
   \end{array}
\label{condition_CFL}
\end{equation}

\noindent 
The boundary condition is supposed to be on the edge of the mesh. The Dirichlet 
boundary condition is written as~: $ \quad  \displaystyle  {\bf{u}}.{\bf{n}} = 0 \, . $

\bigskip   \noindent {\bf \large 6-3 \quad Numerical acoustic absorbing layers }

\noindent 
We deduce from section 4 
the set of equations to solve in the absorbing
layers~:

\[
\left \lbrace
  \begin{array}{l}
        \displaystyle{ \frac{\partial p_x}{\partial t} + \sigma^* (x) \, p_x 
+ \frac{\partial \xi}{\partial x} = 0 } \\
        \displaystyle{ \frac{\partial p_y}{\partial t} + \sigma^* (y) \, p_y 
+ \frac{\partial \zeta}{\partial y} = 0 } \\
        \displaystyle{ \frac{\partial \xi}{\partial t} + \sigma^* (x) \, \xi 
+ \frac{\partial}{\partial x} (p_x + p_y) = 0 } \\
        \displaystyle{ \frac{\partial \zeta}{\partial t} + \sigma^* (y) \, \zeta 
+ \frac{\partial}{\partial y} (p_x + p_y) = 0 }~,
  \end{array}
\right.
\]

\noindent 
where $\sigma^* (x)$ and $\sigma^* (y)$ are the smoothing functions strictly 
positive inside the absorbing layers. The set of equations is ended by a 
Dirichlet condition on the edge of the whole studied domain~: 
$ \quad  \displaystyle  {\bf{u}}.{\bf{n}} = 0 \, , \quad  $  
where ${\bf{n}}$ is the external normal to the domain. The discretization 
of such a boundary condition for the whole studied domain is~:

\[
 \left \lbrace
   \begin{array}{l}
	\xi_{0,j+\frac{1}{2}}^{n+\frac{1}{2}} = 0, \hspace{0.4cm} 0 \leq 
j \leq J-1, \hspace{0.4cm} n \geq 0 \\
	\xi_{J,j+\frac{1}{2}}^{n+\frac{1}{2}} = 0, \hspace{0.4cm} 0 \leq 
j \leq J-1, \hspace{0.4cm} n \geq 0 \\
	\zeta_{i+\frac{1}{2},0}^{n+\frac{1}{2}} = 0, \hspace{0.4cm} 0 \leq 
i \leq J-1, \hspace{0.4cm} n \geq 0 \\
	\zeta_{i+\frac{1}{2},J}^{n+\frac{1}{2}} = 0, \hspace{0.4cm} 0 \leq 
i \leq J-1, \hspace{0.4cm} n \geq 0 
   \end{array}
 \right.
\]

\smallskip \noindent 
We propose to use the same discretisation as before. In the absorbing layers, 
we have to know $p\,$, $\xi$ and $\zeta$ respectively at times $n+\frac{1}{2}$, 
$n+1$ et $n+1$. We decide to center those values in time, we write~:

\[ \left \lbrace
   \begin{array}{lcl}
        \displaystyle{ p_{i+ \frac{1}{2},j+ \frac{1}{2}}^{n+\frac{1}{2}} } & = &
 \displaystyle{ \frac{1}{2} \left( p_{i+\frac{1}{2},j+ \frac{1}{2}}^{n+1} + 
p_{i+\frac{1}{2},j+ \frac{1}{2}}^{n} \right) } \\ \vspace{1 mm}
        \displaystyle{ \xi_{i,j+ \frac{1}{2}}^{n+1} } & = & \displaystyle{ 
\frac{1}{2} \left( \xi_{i,j+ \frac{1}{2}}^{n+\frac{3}{2}} + \xi_{i,j+ 
\frac{1}{2}}^{n+\frac{1}{2}} \right) } \\ \vspace{1 mm}
        \displaystyle{ \zeta_{i+ \frac{1}{2},j}^{n+1} } & = & \displaystyle{
 \frac{1}{2} \left( \zeta_{i+ \frac{1}{2},j}^{n+\frac{3}{2}} + 
\zeta_{i+ \frac{1}{2},j}^{n+\frac{1}{2}} \right) }~.
   \end{array}
   \right.
\]

\smallskip \noindent 
The numerical scheme, while noting $\displaystyle{ \sigma =
\frac{{\Delta t}}{{\Delta x}} = \frac{{\Delta t}}{{\Delta y}} }$ 
(isotrop meshing), can be written as~:

\[
\left \lbrace
   \begin{array}{l}
        \displaystyle{ {p_x}_{i+ \frac{1}{2},j+ \frac{1}{2}}^{n+1} =  
\frac{2 - \sigma^*_{x} (i+ \frac{1}{2})~{\Delta t}}{2 + \sigma^*_{x} 
(i+ \frac{1}{2})~{\Delta t}}~{p_x}_{i+ \frac{1}{2},j+ \frac{1}{2}}^{n} - 
\frac{2 \sigma \left( \xi_{i+1,j+ \frac{1}{2}}^{n+ \frac{1}{2}} - 
\xi_{i,j+ \frac{1}{2}}^{n+ \frac{1}{2}} \right)}{2 + \sigma^*_{x} 
(i+ \frac{1}{2})~{\Delta t}} } \\
        \displaystyle{ {p_y}_{i+ \frac{1}{2},j+ \frac{1}{2}}^{n+1} = 
\frac{2 - \sigma^*_{y} (j+ \frac{1}{2})~{\Delta t}}{2 + \sigma^*_{y} (j+ 
\frac{1}{2})~{\Delta t}}~{p_y}_{i+ \frac{1}{2},j+ \frac{1}{2}}^{n} - 
\frac{2 \sigma \left( \zeta_{i+ \frac{1}{2},j+1}^{n+ \frac{1}{2}} - 
\zeta_{i+ \frac{1}{2},j}^{n+ \frac{1}{2}} \right)}{2 + \sigma^*_{y} (j+ 
\frac{1}{2})~{\Delta t}} } \\
        \displaystyle{ \xi_{i,j+\frac{1}{2}}^{n+ \frac{3}{2}} = \frac{2 - 
\sigma^*_{x} (i)~{\Delta t}}{2 + \sigma^*_{x} (i)~{\Delta t}}~\xi_{i,j+
\frac{1}{2}}^{n+ \frac{1}{2}} - \frac{2 \sigma}{2 + 
\sigma^*_{x} (i)~{\Delta t}}~\left( p_{i+\frac{1}{2},j+\frac{1}{2}}^{n+1} - 
p_{i-\frac{1}{2},j+\frac{1}{2}}^{n+1} \right) } \\
        \displaystyle{ \zeta_{i+\frac{1}{2},j}^{n+ \frac{3}{2}} = 
\frac{2 - \sigma^*_{y} (j)~{\Delta t}}{2 + 
\sigma^*_{y} (j)~{\Delta t}}~\zeta_{i+\frac{1}{2},j}^{n+ \frac{1}{2}} -
 \frac{2 \sigma}{2 + \sigma^*_{y} (j)~{\Delta t}}~\left( p_{i+\frac{1}{2},
j+\frac{1}{2}}^{n+1} - p_{i+\frac{1}{2},j-\frac{1}{2}}^{n+1} \right) }~,
    \end{array}
\right.
\]

\noindent with $p \equiv p_x+p_y$~.

\bigskip \bigskip  \newpage \noindent {\bf \large 7) \quad Numerical tests}  

\noindent 
The computational domain is shared into two areas~: the studying medium and 
the absorbing layers. Those two areas are defined as shown below~: 

\centerline { 
  \begin{picture}(100,68)
        \put(4,4){\line(0,1){60}}
        \put(4,4){\line(1,0){60}}
        \put(64,64){\line(0,-1){60}}
        \put(64,64){\line(-1,0){60}}
        \put(22,22){\line(0,1){24}}
        \put(22,22){\line(1,0){24}}
        \put(46,46){\line(0,-1){24}}
        \put(46,46){\line(-1,0){24}}
        \put(33,34){\line(1,0){2}}
        \put(34,33){\line(0,1){2}}
        \put(79,59){\vector(-1,0){25}}
        \put(80,58){Absorbing layers}
        \put(79,26){\vector(-1,0){40}}
        \put(80,25){Studying medium}
        \put(34,3){{\line(0,1){1}}}
        \put(33,0){0}
        \put(46,3){{\line(0,1){1}}}
        \put(44,0){$x_{max}$}
        \put(64,3){{\line(0,1){1}}}
        \put(64,0){$x_{max} + L_{pml}$}
        \put(22,3){{\line(0,1){1}}}
        \put(16,0){$-x_{max}$}
        \put(3,46){{\line(1,0){1}}}
        \put(-6,46){$x_{max}$}
        \put(3,22){{\line(1,0){1}}}
        \put(-9,22){$-x_{max}$}        
  \end{picture} }

\bigskip   \noindent {\bf \large 7-1 \quad Mathematical experiments}

\noindent 
We first propose to numerically validate our mathematical analysis of the absorbing 
layers acoustic system (\ref{pml_matrice}) with a test case proposed by O. Pironneau 
\cite{pironneau}. We decide to place an acoustic pulse {\bf{inside}} the absorbing 
layers, in the direction of the eigenvector associated to the eigenvalue 0 in order 
to enforce the unstability due to the lack of hyperbolicity as shown before. For this
test, the computational domain, symetric in the $x$-direction and the $y$-direction, 
is defined by $x_{max} = 5$ and $L_{pml} = 45$. We consider two different tests. 
For the first one, we solve~:

\moneq
\left \lbrace
    \begin{array}{l}
        \displaystyle{ \frac{\partial W}{\partial t} + A \frac{\partial W}{\partial x} 
+ B \frac{\partial W}{\partial y} + C \, W = \psi(t) \, \delta_{x_a,y_a} \, V_{M^2} } \\
        W (0) = 0~.
    \end{array}
\right.
\label{PbM1}
\end{equation}

\noindent For the second one, we solve~:

\moneq
\left \lbrace
    \begin{array}{l}
        \displaystyle{ \frac{\partial W}{\partial t} + A 
\frac{\partial W}{\partial x} + B \frac{\partial W}{\partial y} + C \, W = 0~, } \\
        \displaystyle{ W (0) = 
\left(
  \begin{array}{c}
        0 \\ 0 \\ \exp \left( - (\ln 2) \frac{(x-x_a)^2 + (y-y_a)^2}{9} \right) \\ 
\exp \left( - (\ln 2) \frac{(x-x_a)^2 + (y-y_a)^2}{9} \right)
  \end{array}
\right)~. }
    \end{array}
\right.
\label{PbM2}
\end{equation}

\smallskip \noindent 
For the first test case, the excitation $\psi(t) \, V_{M^2}$ is centered {\bf{inside}} 
the absorbing layers at $(x_a,y_a) = (25,0)$. The absorbing coefficients in the absorbing 
layers are constant in the $x$ and in the $y$-direction. The reference excitation is~:

\smallskip \noindent \centerline {$ \displaystyle 
        \psi_{xy} = \exp \left( - (\ln 2) \frac{(x-x_a)^2 + (y-y_a)^2}{9} \right)~. \, $}

\smallskip \noindent 
We solve the problem (\ref{PbM1}) with the particular data given by~:

\[
  \left \lbrace
    \begin{array}{l}
        \check{p} (x_a,y_a,t) = 0 \\
        \check{\xi} (x_a,y_a,t) = \frac{\partial \psi_{xy}}{\partial y} \\
        \check{\zeta} (x_a,y_a,t) = -\frac{\partial \psi_{xy}}{\partial x}~.
    \end{array}
  \right.
\]
The observing points are centered at $(x_1,y_1) = (45,0)$, $(x_2,y_2) = 
(25,0)$, $(x_3,y_3) = (0,0)$, $(x_4,y_4) = (-45,0)$, $(x_5,y_5) = (0,25)$ and $(x_6,y_6) 
= (0,-25)$. We observe the results at $(x_2,y_2) = (25,0)$. We obtain the graph presented
of Figure 1. 

\smallskip  
\centerline { \includegraphics{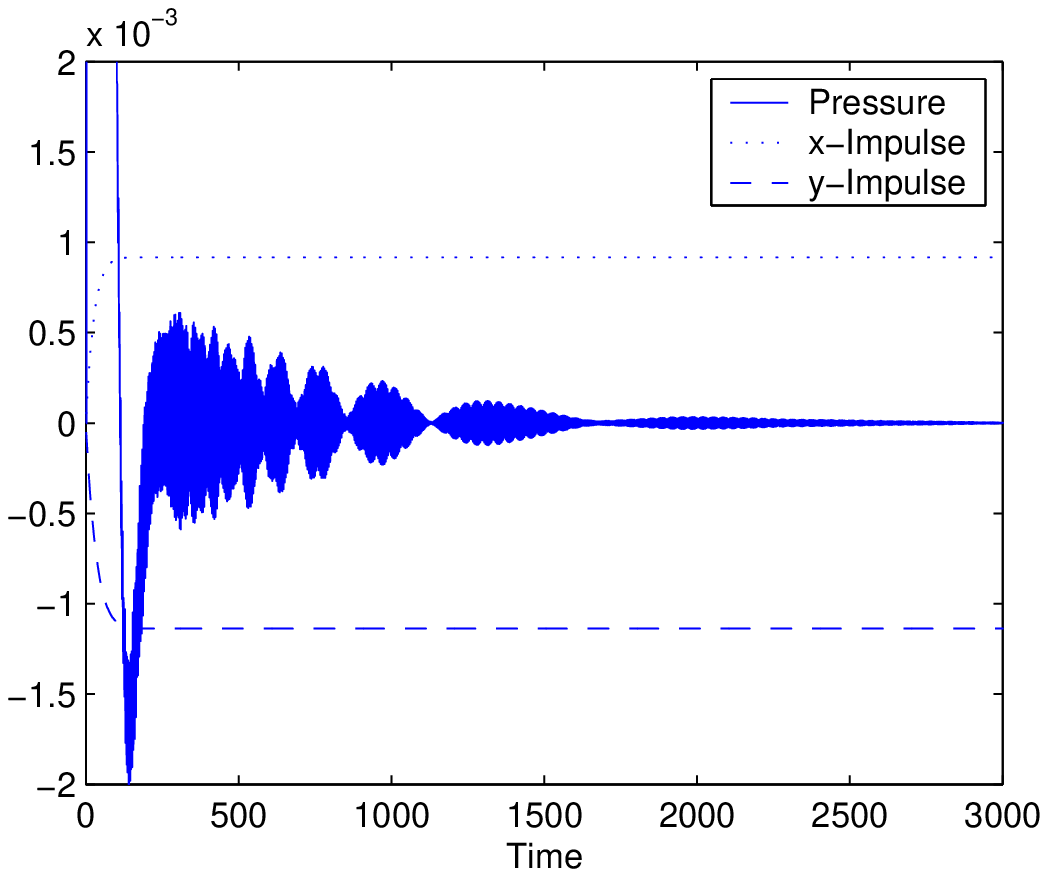} }
\smallskip  \centerline { {\bf Figure 1}. \quad  
 Pressure and impulse fields at $(x,y) = (25,0)$. }
\bigskip 

\smallskip \noindent 
We remark that, for long times, the pressure field converges to 0 and the impulse 
fields to a non-zero value. We obtain the results predicted by our simple 1-D model 
(see proposition \ref{prop_1D}) noting $u$ the pressure field and $v$ the impulse 
fields $\xi$ or $\zeta$. We notice that the sign is changing if we consider a positive 
or a negative source. We observe approximately the same results for the other observing 
points. The pressure field converges to 0 except at $(x_5,y_5)$ and $(x_6,y_6)$, 
where there is a slight residual $r_p$ at $(x_5,y_5)$ and $-r_p$ at $(x_6,y_6)$. 
The $\xi$ impulse field converges to 0 at $(x_1,y_1)$, $(x_3,y_3)$ and $(x_4,y_4)$. 
The residual at $(x_5,y_5)$ is $r_{\xi}$ and $-r_{\xi}$ at $(x_6,y_6)$. The $\zeta$ 
impulse field always converges to a constant value as predicted by the 1-D model. \\

\smallskip \noindent 
For the second test, we solve the problem (\ref{PbM2}) with an excitation centered 
{\bf{inside}} the absorbing layers at $(x_a,y_a) = (25,0)$. As before, we analyse 
the results at $(x_2,y_2) = (25,0)$. The evolution of the pressure and impulse fields 
are drawn in the following figure~:

\bigskip 
\centerline {
{\includegraphics [width=7cm] {{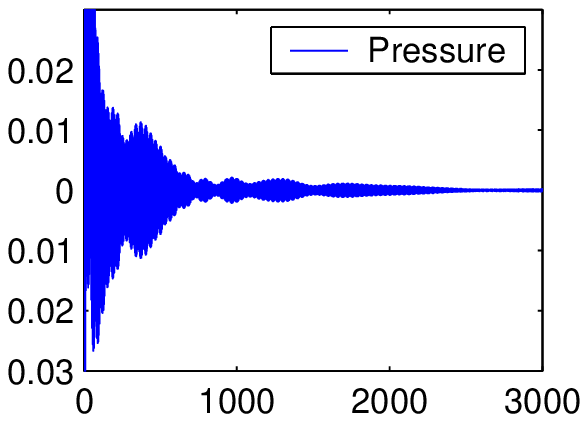}}}
{\includegraphics [width=7cm]{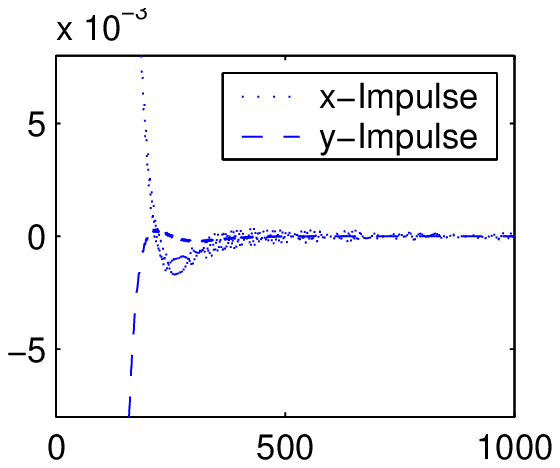}}} 
\smallskip  \centerline { {\bf Figure 2}. \quad  
 Pressure and impulse fields at $(x,y) = (25,0)$. }
\bigskip 

\smallskip \noindent 
All the fields converge to 0 for all the observing points. We obtain the same results 
as predicted by the proposition \ref{prop_1D} with an excitation function $\psi (t)$ 
such as $\int_0^\infty \psi (t) {\rm d} t$ is bounded. \\

\smallskip \noindent 
The conclusion of these mathematical experiments is that even if we excite this 
non hyperbolic system in the direction of the non caracteristic vector, the zero 
order damping terms insure that there is no numerical explosion of our results.

\bigskip  \noindent {\bf \large 7-2 \quad Physical experiments} 

\noindent 
We have proven the stability of our absorbing layers. We now want to study numerical 
reflections of outgoing waves from the boundaries of the computational domain for 
various speeds of the external flow. We then consider the problem (\ref{sys2_2}) in 
the computational domain defined by $x_{max} = 25$ with the absorbing layers outside 
($L_{pml}$ is now a parameter), an acoustic pulse centered at $(x_a,y_a) = (0,0)$ and 
an excitation in the right hand side $(\check{p},\check{\xi},\check{\zeta})^t$ given by~:

\vspace*{-5pt}
\moneq
\left \lbrace
  \begin{array}{l}
	\displaystyle{ \check{p} \, (x,y,t) = \exp \left( - (\ln 2) 
\frac{(x-x_a)^2 + (y-y_a)^2}{9} \right) \sin (\pi t) } \\
	\check{\xi} (x,y,t) = 0 \\
	\check{\zeta} (x,y,t) = 0~. 
  \end{array}
\right.
\label{excitation}
\end{equation}

\smallskip  \noindent 
We compare the calculated solution, denoted by $p$, to the numerical solution obtained 
in the domain defined by $x_{max} = 150$. We take ${\Delta x} = 1$ and ${\Delta t}$ 
following the CFL (\ref{condition_CFL}). For $t < 300 \, {\Delta t}$, it is easy to see 
that no reflection from the boundaries can interact with the solution within the small 
domain $[-25,25]^2$ and such a solution is the numerical solution in an infinite domain 
for $t < 300 \, {\Delta t}$. This solution is considered as a reference, noted $p_{ref}$, 
and the computation of the error $|p - p_{ref}|$ for each time step indicates the
efficiency 
of the absorbing layers. As explained in Section 5, 
the boundary condition is 
imposed on the edge of a cell and the pressure is calculated in the middle of the cell. 
The observing point is then taken only half a cell near the absorbing layers at $(25,0)$.

\bigskip    \smallskip  
\centerline {	\includegraphics  [width=10.0cm]    {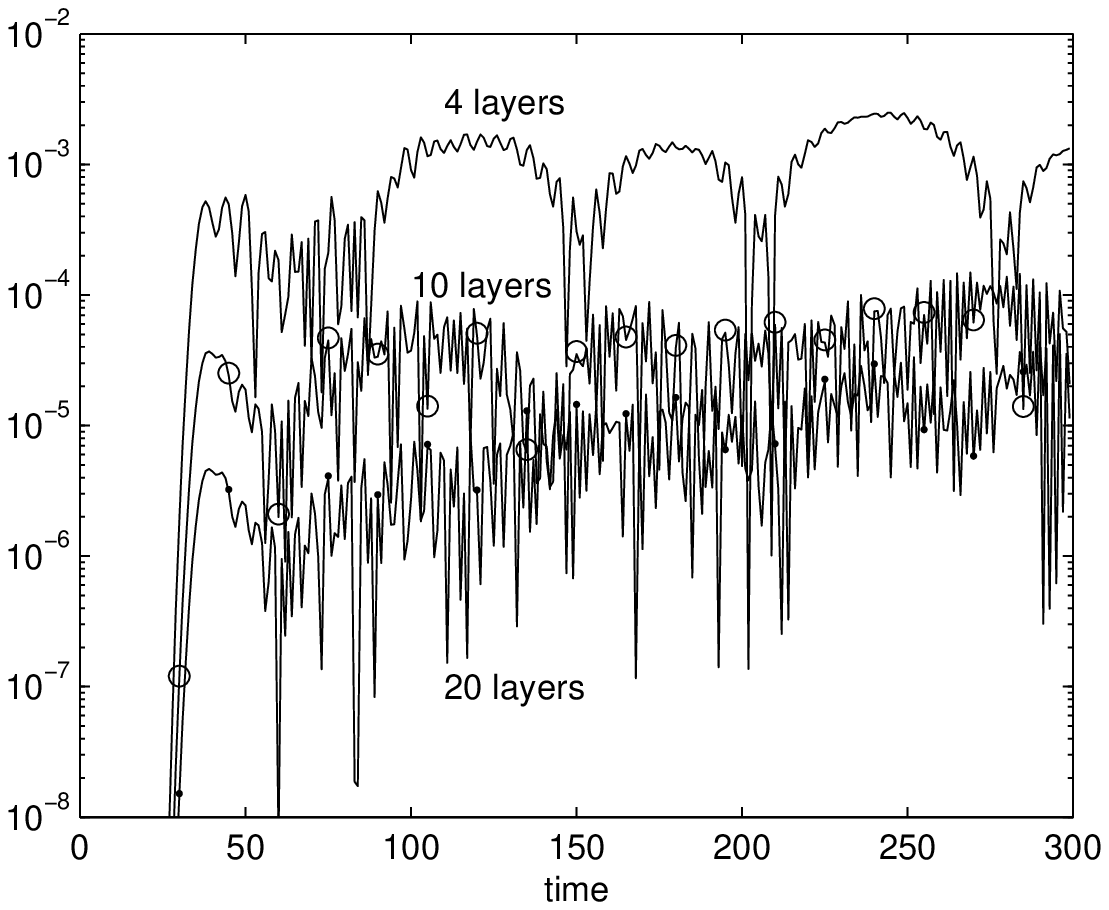}}

 \noindent  {\bf Figure 3}. \quad 
 ${\rm L}_2-$error  of the pressure for 4, 10 and 20 absorbing 
layers, $\frac{{\bf{u}}}{c_0} = (0.5, 0)$. 
\bigskip 

 \smallskip  
\centerline { 	\includegraphics  [width=10.0cm]  {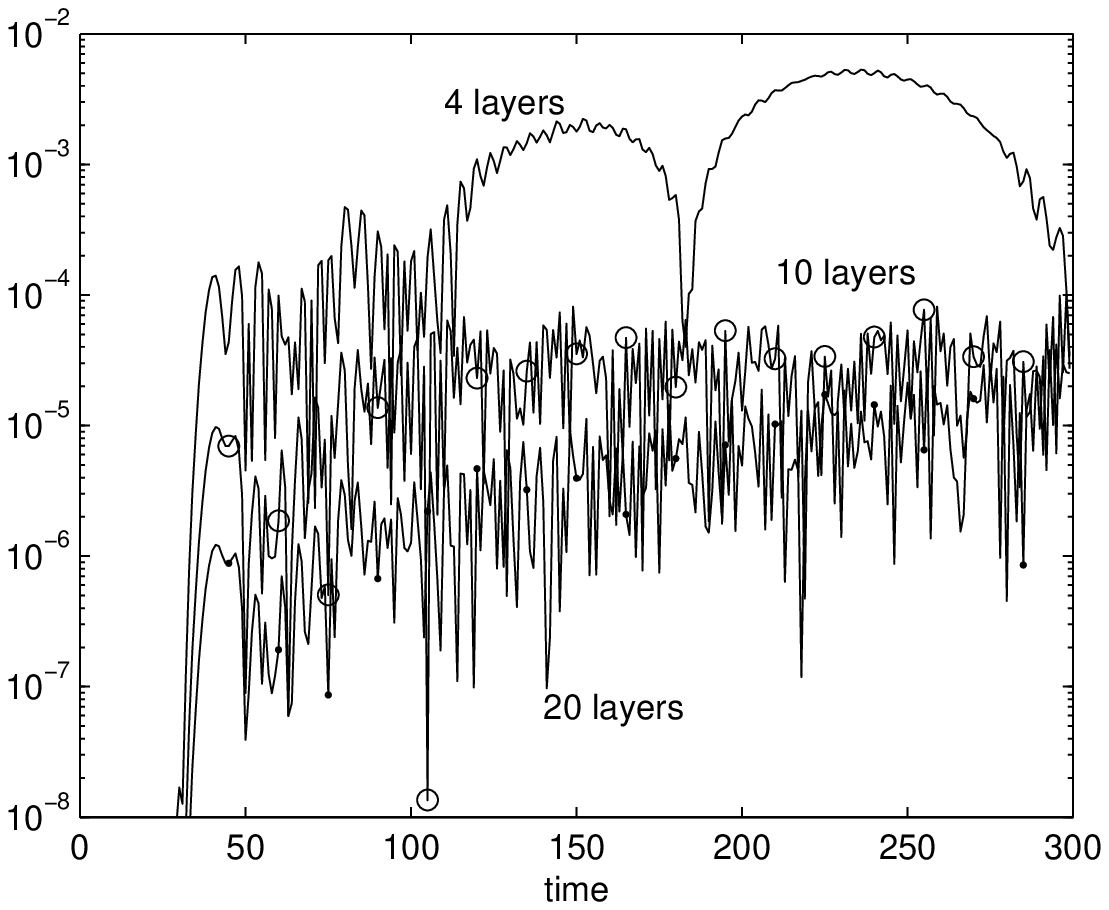}}

  \noindent  {\bf Figure 4}. \quad  
 ${\rm L}_2-$error  of the pressure for 4, 10 and 20 absorbing 
layers, $\frac{{\bf{u}}}{c_0} = (0.5, 0)$. 

\newpage  

~   

\bigskip  \bigskip  \bigskip  
  \smallskip    
\centerline {	\includegraphics  [width=10.0cm]    {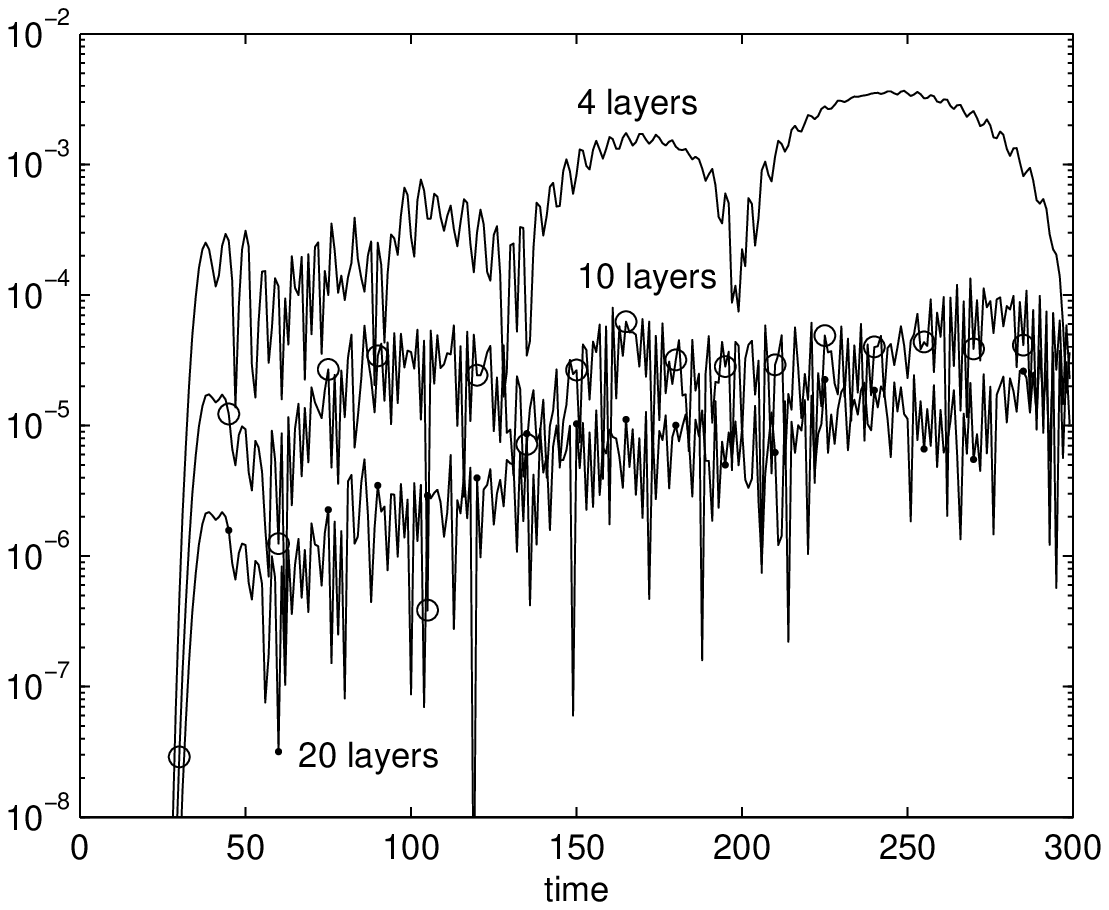}}

 \noindent  {\bf Figure 5}. \quad  
 ${\rm L}_2-$error of the pressure for 4, 10 and 20 absorbing layers, 
$\frac{{\bf{u}}}{c_0} = (\frac{1}{2 \sqrt{2}} , \frac{1}{2 \sqrt{2}})$.

\smallskip  
\centerline {	\includegraphics [width=10.0cm]   {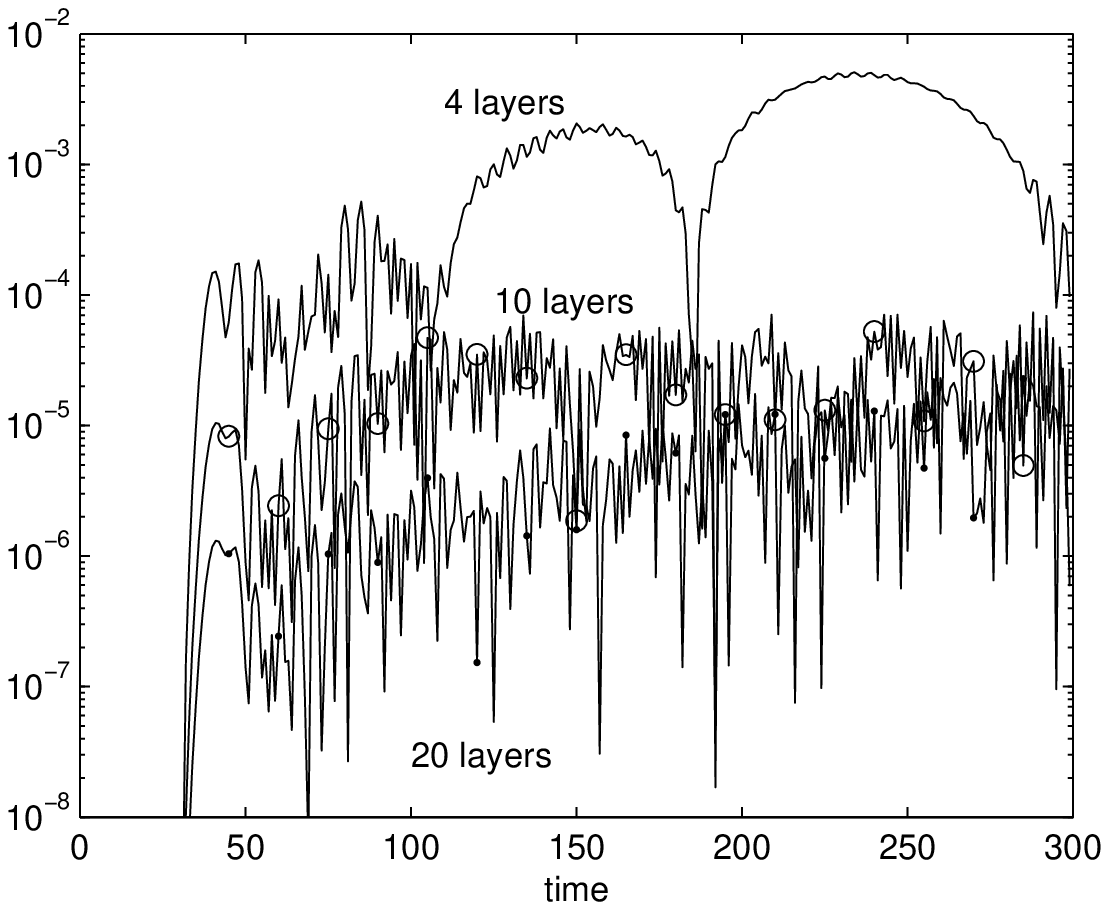} }

 \noindent  {\bf Figure 6}. \quad  
  ${\rm L}_2-$error of the pressure for 4, 10 and 20 absorbing layers, 
$\frac{{\bf{u}}}{c_0} = (\frac{2}{\sqrt{17}} , \frac{1}{2 \sqrt{17}})$.
 
 \newpage   \centerline {
{\includegraphics [width=6.9cm] {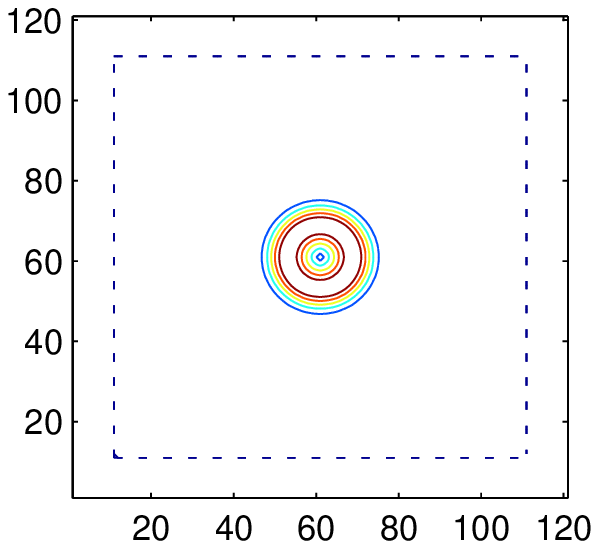}}   
{\includegraphics [width=6.9cm] {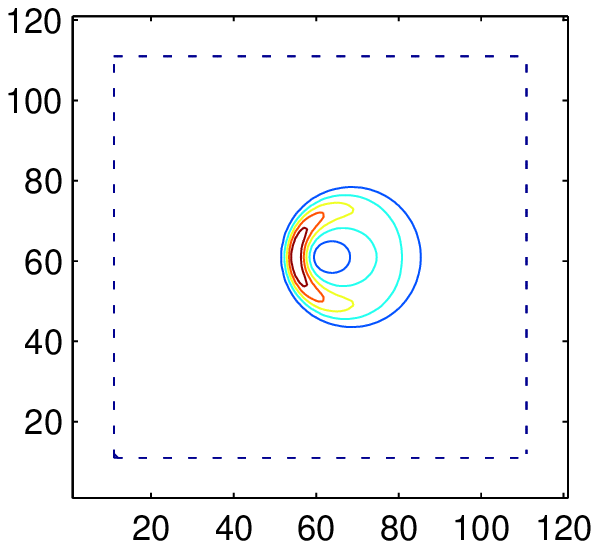}}}
  \bigskip  \centerline {
{\includegraphics [width=6.9cm] {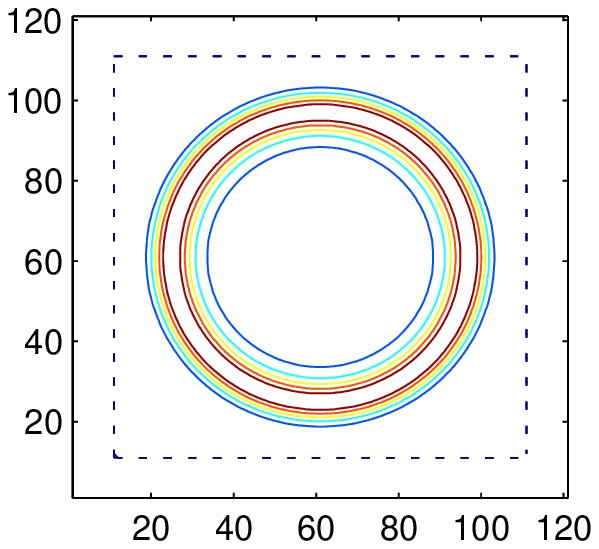}}  
{\includegraphics [width=6.9cm] {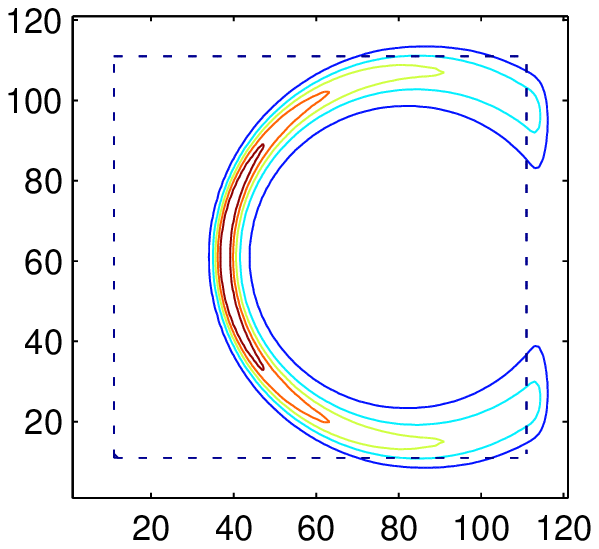}}} 
 \bigskip  \centerline { 
{\includegraphics [width=6.9cm] {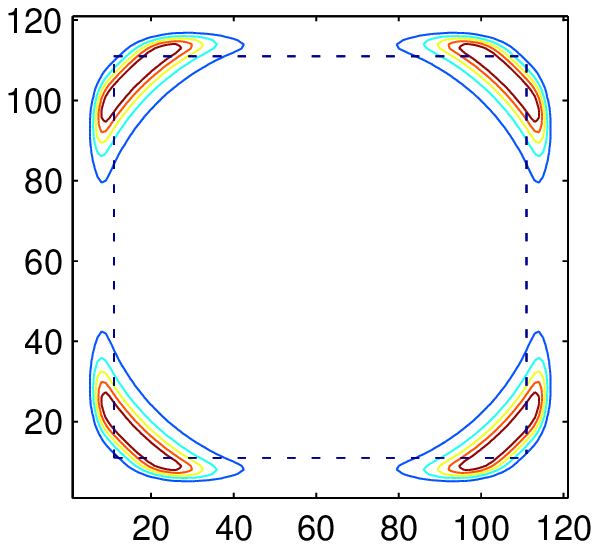}}   
{\includegraphics [width=6.9cm] {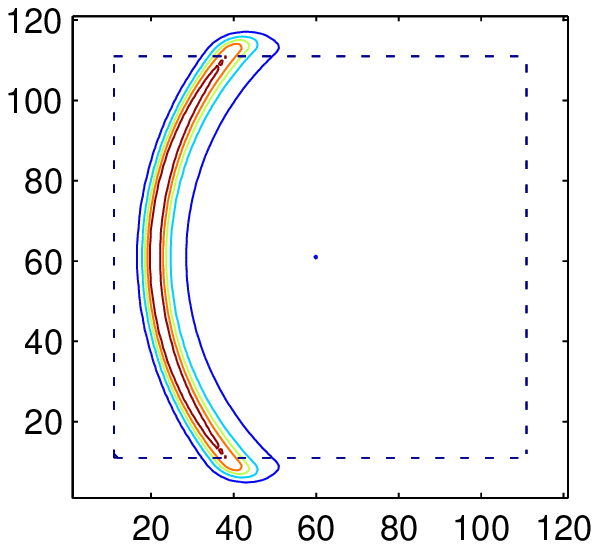}}}
\bigskip  \noindent  {\bf Figure 7}. \quad   
Iso-lines of the pressure field 
for $t=40\,{\Delta t}$, $t=80\,{\Delta t}$ and $t=120\,{\Delta t}$.

 \newpage   \centerline { 
{\includegraphics [width=6.9cm] {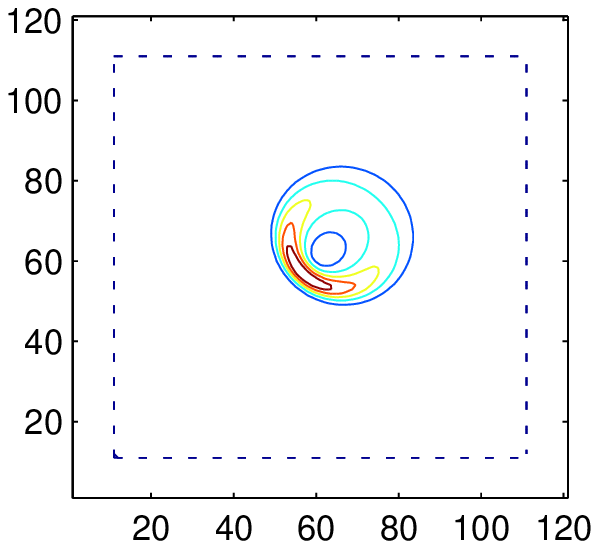}} 
{\includegraphics [width=6.9cm] {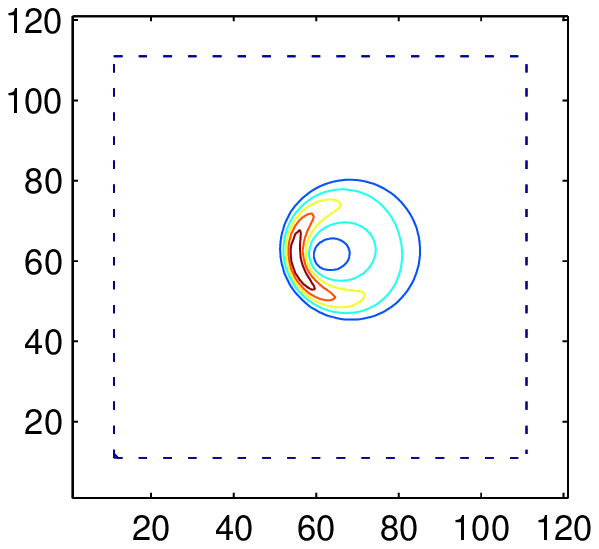}}}
  \bigskip  \centerline { 
{\includegraphics [width=6.9cm] {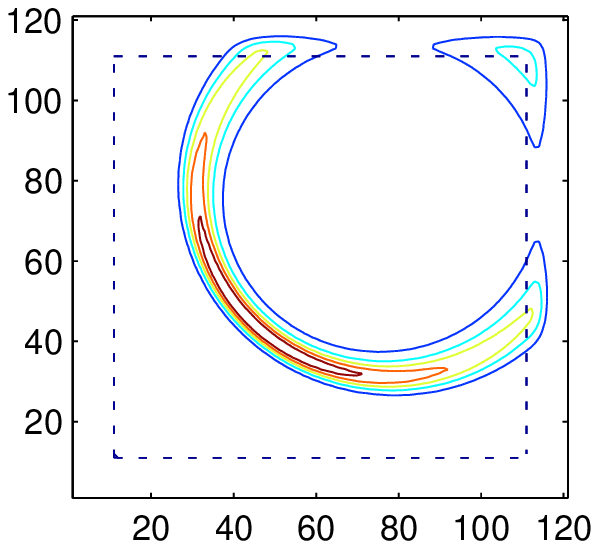}}   
{\includegraphics [width=6.9cm] {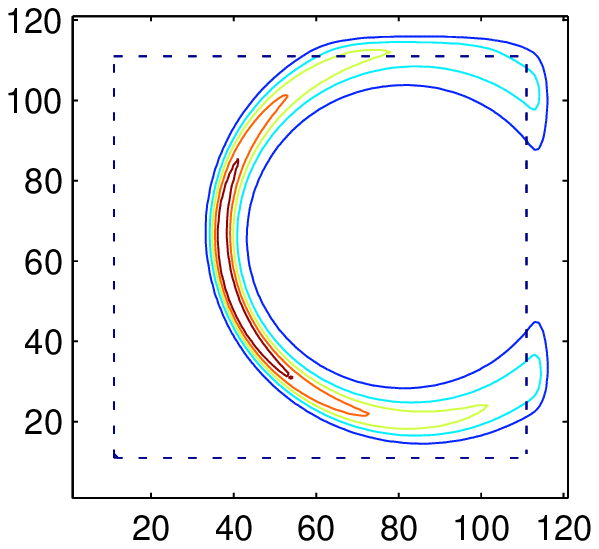}} } 
  \bigskip  \centerline { 
{\includegraphics [width=6.9cm] {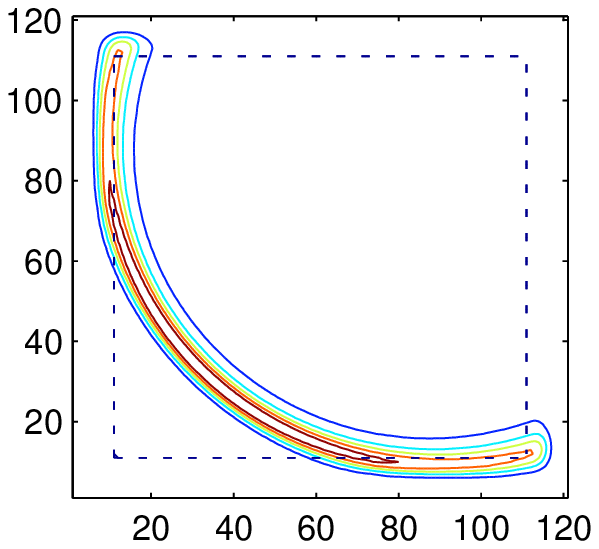}}   
{\includegraphics [width=6.9cm] {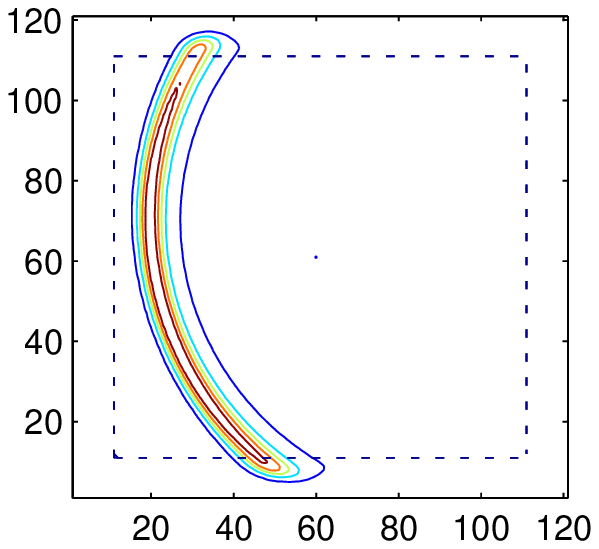}} }
\bigskip  \noindent  {\bf Figure 8}. \quad  Iso-lines of the pressure field 
for $t=40\,{\Delta t}$, $t=80\,{\Delta t}$ and $t=120\,{\Delta t}$. 

 \newpage

\smallskip \noindent 
We first consider a domain without external flow. We notice that the $L_2$-error of the 
pressure field $p$ computed for each time step at the observing point $(25,0)$ for various 
thickness for the absorbing layers. We notice that for small absorbing layers (4 cells), 
we have a good accuracy for the results, and that increasing the thickness from 4 cells 
to 20 cells improves the accuracy by 2 orders of magnitude.

\smallskip \noindent 
We also compare, for various thickness of the absorbing layers, the ``exact'' and the
numerical solution considering the same observing point and the same acoustic source, for
three velocity vectors defined by ${\bf{u}} = (0.5 \, c_0, 0)$, ${\bf{u}} = (\frac{1}{2
  \sqrt{2}} \, c_0, \frac{1}{2 \sqrt 2} \, c_0)$ and ${\bf{u}} = (\frac{2}{\sqrt{17}} \,
c_0, \frac{1}{2 \sqrt{17}} \, c_0)$. We notice an improvement of the accuracy by 2 orders
of magnitude for an absorbing layers growing from 4 cells to a 20 cells.

\smallskip \noindent 
The results are satisfying. Nevertheless, when the number of cells in the absorbing 
layers is increasing, the error is small but remains measurable, even for the long 
times. We think that this behavior could be improved in future work.

\bigskip \bigskip  \noindent {\bf \large  Conclusion }   

\noindent 
We have explored a new method for solving the equations of advective acoustics based
 on a change a space-time variables (Lorentz transform) and a change of unknown 
variables. We have also derived a system of equations (\ref{sys_pml_uv}) to modelize 
the absorbing layers for the acoustic model. The system of partial differential 
equations established in the absorbing layers is well-posed due to the zero order term.
 The staggered grid ``HaWAY'' method has been used for the numerical implementation and 
experiments have proven the efficiency of such a method. When we force a punctual 
acoustic source inside this numerical domain, our experiments show that the results 
remain bounded and our method is stable from a practical point of view. Notice that 
we explain our method only in two-dimensional space, but we extend it easily to 
three-dimensional space.

\bigskip \bigskip   \noindent {\bf \large Acknowledgments}   

\noindent   
The authors thank L. Halpern and A. Rahmouni for stimulating scientific 
discussions, O. Pironneau for suggesting an interesting mathematical test 
and European Aeronautics Defence and Space at Suresnes and Airbus France 
for financial support of this reseach.

\bigskip \bigskip  

\end{document}